# SBP-SAT finite difference discretization of acoustic wave equations on staggered block-wise uniform grids

Longfei Gao[a,*], David C. Del Rey Fernández[b,c], Mark Carpenter[b], David Keyes[a]

[a]Division of Computer, Electrical and Mathematical Sciences and Engineering, King Abdullah University of Science and Technology, Thuwal 23955-6900, Saudi Arabia.
[b]NASA Langley Research Center, Hampton, VA 23681, USA
[c]National Institute of Aerospace, Hampton, VA 23666, USA

**Abstract**

We consider the numerical simulation of the acoustic wave equations arising from seismic applications, for which staggered grid finite difference methods are popular choices due to their simplicity and efficiency. We relax the uniform grid restriction on finite difference methods and allow the grids to be block-wise uniform with nonconforming interfaces. In doing so, variations in the wave speeds of the subterranean media can be accounted for more efficiently. Staggered grid finite difference operators satisfying the summation-by-parts (SBP) property are devised to approximate the spatial derivatives appearing in the acoustic wave equation. These operators are applied within each block independently. The coupling between blocks is achieved through simultaneous approximation terms (SATs), which impose the interface condition weakly, i.e., by penalty. Ratio of the grid spacing of neighboring blocks is allowed to be rational number, for which specially designed interpolation formulas are presented. These interpolation formulas constitute key pieces of the simultaneous approximation terms. The overall discretization is shown to be energy-conserving and examined on test cases of both theoretical and practical interests, delivering accurate and stable simulation results.

*Keywords:* Summation by parts, simultaneous approximation terms, seismic wave modeling, nonconforming interface, staggered grid, long time instability

## 1. Introduction

Numerical simulation of wave phenomena, commonly referred to as forward modeling in geophysics, is of vital importance in seismic applications. In modern-day industrial practice, computationally synthesized solutions are compared against the recorded field data to drive the iterative inversion processes that characterize the subsurface structures. In most modern seismic inversion approaches, cf. [1–4], independent forward modeling tasks need to be performed for various source terms at each inversion iteration step, which constitute the most computationally intensive part of the overall process.

The forward modeling problem can be posed in either the time domain or the frequency domain. Here, we focus on the time domain approach due to its robustness and efficiency for a wide range of applications. Various discretization methods have been proposed for this task, including finite element methods on unstructured grids, e.g., [5–7], and finite difference methods on uniform grids, e.g., [8–10]. Although finite element methods are more flexible for handling complicated geometry and local refinement, finite difference methods remain widely employed in the seismic community due to their simplicity and efficiency.

However, efficiency of finite difference methods can be impaired due to large variations in the model parameters. Taking the acoustic wave equation as an example, the spatial sampling rate of a uniform grid is determined by the accuracy requirement, which is restricted by the minimum wave speed throughout the overall medium. On the other hand, the temporal sampling rate is restricted by the Courant-Friedrichs-Lewy stability condition, which is determined by the spatial sampling rate and the maximum wave speed throughout the overall medium. Consequently, variations in the wave speed lead to spatial oversampling in domains of high wave speed and temporal oversampling in domains of low wave speed.

*Corresponding author
*Email addresses:* `longfei.gao@kaust.edu.sa` (Longfei Gao), `dcdelrey@gmail.com` (David C. Del Rey Fernández), `mark.h.carpenter@nasa.gov` (Mark Carpenter), `david.keyes@kaust.edu.sa` (David Keyes)

Motivated by these practical issues, we aim to increase the flexibility of finite difference methods by allowing the underlying grids to be block-wise uniform, connected by nonconforming interfaces, so that domains with different wave speeds can be sampled at their appropriate rates. Due to geological sedimentation and consolidation, wave speeds in subterranean media tend to increase with depth. Therefore, we are particularly interested in grids that have larger grid spacing for the deeper part of the simulation domain. Earlier attempts exist in the geophysical community, e.g., [11–13]. Unfortunately, the discretization methods therein demonstrate unstable behavior that manifests for long time simulation.

In this work, we propose finite difference discretization of the acoustic wave equation on staggered block-wise uniform grid based on summation-by-parts (SBP) operators and simultaneous approximation terms (SATs). The SBP operators date back to [14], and since then have been continuously enriched in both theory and application for various physical problems. The SATs technique was originally proposed in [15] as a tool to weakly impose the boundary condition and has since then become the standard boundary or interface treatment approach to combine with the SBP operators. Readers may consult the two review papers [16, 17] and the references therein for more details about the SBP-SAT framework.

By discrete energy analysis, the finite difference discretization proposed in this work is shown to be energy-conserving, thus eliminating the unstable behavior encountered in previous attempts. Moreover, in previous attempts, ratio of the grid spacing of neighboring blocks is chosen to be an integer to ease the effort of grid coupling. In this work, with the flexibility offered by the SATs technique and specially designed interpolation operators, we are able to relax this ratio to allow rational numbers.

The rest of this paper is organized as follows. In Section 2, we present the background of the physical problem and the long time instability issue that can appear in nonuniform grid simulations as motivation for this study. In Section 3, we present the SBP operators for staggered grid finite difference discretization and the SATs for interface coupling. Discrete energy analysis is used to demonstrate the energy-conserving property of the proposed discretization method. In Section 4, numerical examples are presented to demonstrate the efficacy of the proposed method. In Section 5, we briefly remark on several relevant issues that are not the focus of this work, but can be of practical interest to readers. References are pointed out for these issues. Finally, we conclude this work with Section 6.

## 2. Problem Description

We consider the 2D acoustic wave equation posed as the following first-order system:

$$\left\{\begin{array}{rcl} \dfrac{\partial p}{\partial t} & = & -\rho c^2\left(\dfrac{\partial u}{\partial x}+\dfrac{\partial v}{\partial y}\right)+\mathcal{S}; \\ \dfrac{\partial u}{\partial t} & = & -\dfrac{1}{\rho}\dfrac{\partial p}{\partial x}; \\ \dfrac{\partial v}{\partial t} & = & -\dfrac{1}{\rho}\dfrac{\partial p}{\partial y}, \end{array}\right. \tag{1}$$

where $p$, $u$ and $v$, standing for pressure, particle velocity in $x$-direction and particle velocity in $y$-direction, respectively, are the sought solution variables. Parameters of this equation include density $\rho$ and wave speed $c$. The source term that drives the wave propagation is denoted by $\mathcal{S}$. In the upcoming analysis, the homogeneous version of (1) is often considered, where the source term $\mathcal{S}$ is omitted. The medium is assumed to be at rest at the beginning of the simulation. Translating into initial conditions, this means that $p$, $u$ and $v$, as well as their derivatives, are zero.

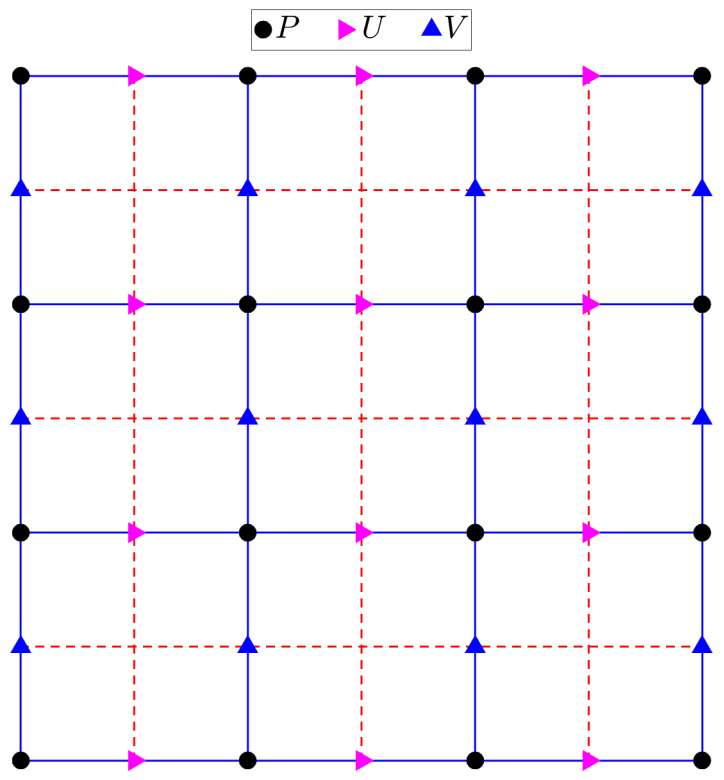

Figure 1: Illustration of the 2D uniform staggered grid.

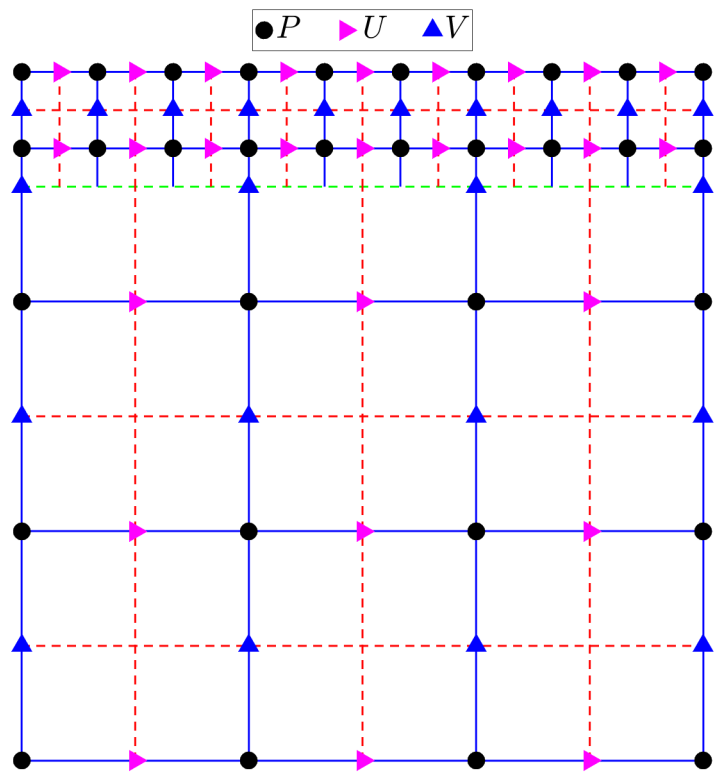

Figure 2: Illustration of the 2D block-wise uniform staggered grid with a nonconforming interface.

In the seismic community, system (1) is often solved on the staggered grid, as demonstrated in Figure 1, where each solution variable occupies its own subgrid. The $u$ subgrid is right shifted from the $p$ subgrid for half the grid spacing. For instance, the $i$th grid point of the $u$ subgrid is at the midpoint between the $i$th and the $(i+1)$th grid points of the $p$ subgrid. Similarly, the $v$ subgrid is up shifted from the $p$ subgrid for half the grid spacing.

As mentioned in the introduction, we are interested in the type of grids that have larger grid spacing for deeper part of the simulation domain, as illustrated in Figure 2 for a case of staggered block-wise uniform grid with three times larger grid spacing below the interface. (This 3:1 ratio of grid spacing, or other odd numbers, is favored for staggered grids since the grids of coarser resolution can be considered as a subset of the grids of finer resolution, cf. [12].) Earlier attempts on finite difference simulation of seismic waves on this type of grids exist, e.g., [11, 12]. Unfortunately, these earlier attempts suffer from the so-called long time instability issue, which appears only at the late stage of the simulation and has little visible influence at the early stage.

To illustrate, we solve equation (1) on a block-wise uniform grid as illustrated in Figure 2. Numerical configuration of this test is identical to that described in [18, pp. 1099-1100], except that fourth-order staggered central differences are used in the coarse region to improve accuracy. Time history of pressure is recorded at one point, i.e., the receiver location, and displayed in Figures 3a and 3b for the first and last two seconds, respectively. We observe that at the beginning of the simulation, the numerical result simulated on the block-wise uniform grid matches very well with the uniform grid result. However, at the later stage of the simulation, the highly oscillatory unphysical modes dominate the block-wise uniform grid result. It is reported in [18] that these unphysical modes grow at exponential rates. Various a posteriori numerical techniques have been proposed to mitigate this long time instability issue, e.g., spatial filtering [11, 12] and temporal filtering [19]. In the upcoming sections, we develop discretization operators that are inherently stable and hence prevent the kind of unstable behavior demonstrated in Figure 3b.

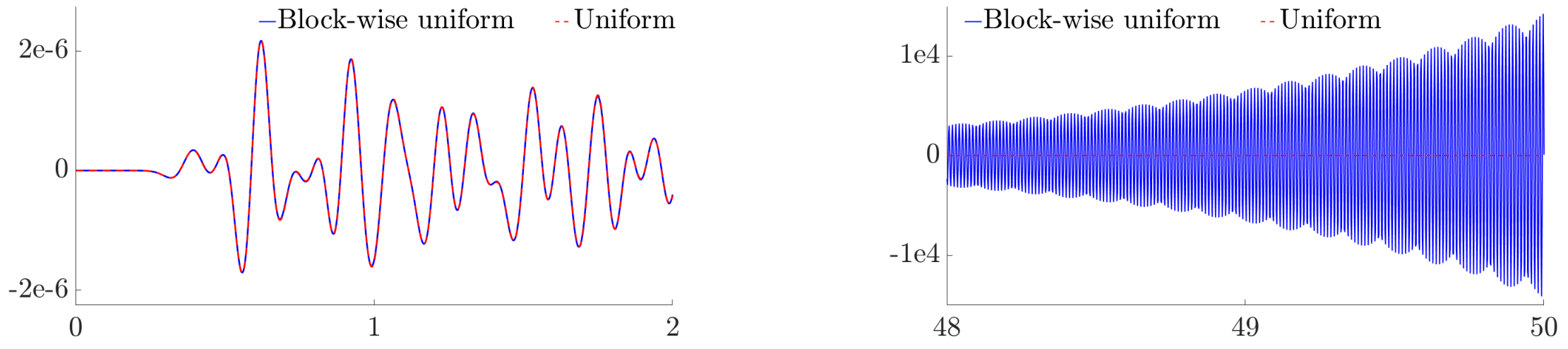

(a) First two seconds.

(b) Last two seconds.

Figure 3: Long time instability on block-wise uniform grid with a homogeneous medium.

## 3. Methodology

In this section, we present the methodology of the SBP-SAT discretization framework and demonstrate how it leads to energy-conserving discretization of the acoustic wave equation on staggered block-wise uniform grids.

### 3.1. A 1D model problem

To start, we consider the following 1D prototype PDE system:

$$\begin{cases} \dfrac{\partial p}{\partial t} = -\dfrac{\partial v}{\partial x}; \\ \dfrac{\partial v}{\partial t} = -\dfrac{\partial p}{\partial x}, \end{cases} \tag{2}$$

defined on interval $(x_L, x_R)$ and discretized on a uniform grid with grid spacing $\Delta x$, as demonstrated in Figure 4, where the $v$ subgrid is staggered to the right of the $p$ subgrid. The particular form of the boundary condition is left unaddressed at this stage. We first present the spatial discretization operators for (2) that satisfy the SBP property in Section 3.1.1.

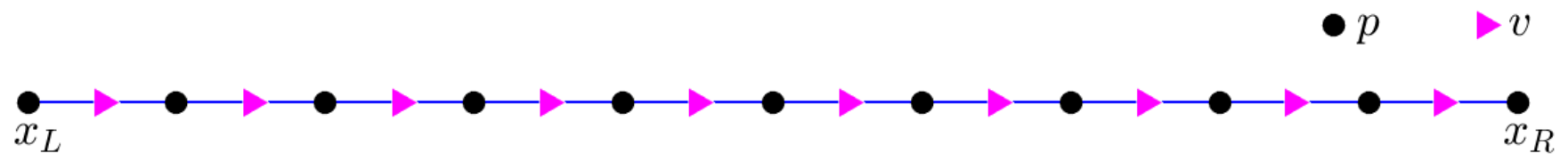

Figure 4: Illustration of the 1D staggered grid.

We note here that in [20], the authors have presented SBP finite difference operators on staggered grids where both subgrids are aligned on the boundaries, i.e., two extra $v$ subgrid points are attached at $x_L$ and $x_R$, comparing to Figure 4. In this work, we do not make this adjustment and the SBP property is retained by invoking the concept of projection in designing the finite difference operators, to be explained in detail in Section 3.1.1. Other usage of staggered grids within the SBP-SAT framework can be found in, e.g., [21, 22].

#### 3.1.1. SBP discretization operators for the 1D model problem

We aim to find spatial discretization operators that mimic the continuous operators $\partial v/\partial x$ and $\partial p/\partial x$ in the context of energy analysis. At the continuous level, we define the energy associated with (2) as

$$\mathscr{E} = \frac{1}{2}\int_{x_L}^{x_R} \left(p^2 + v^2\right) dx. \tag{3}$$

Taking the temporal derivative of $\mathscr{E}$, substituting in equation (2) and then applying integration by parts, we obtain

$$\frac{d\mathscr{E}}{dt} = \int_{x_L}^{x_R} \left(p\frac{\partial p}{\partial t} + v\frac{\partial v}{\partial t}\right) dx = \int_{x_L}^{x_R} -\left(p\frac{\partial v}{\partial x} + v\frac{\partial p}{\partial x}\right) dx = -(pv)\Big|_{x_L}^{x_R}. \tag{4}$$

According to (4), the evolution of the energy associated with system (2), as given in (3), depends only on the solution at the boundary. We aim to mimic this property with specially designed spatial discretization operators, which will be referred to as the SBP operators.

We denote the semi-discretized version of (2) as

$$\begin{cases} \mathcal{A}^P \dfrac{dP}{dt} = -\mathcal{A}^P \mathcal{D}^V V; \\ \mathcal{A}^V \dfrac{dV}{dt} = -\mathcal{A}^V \mathcal{D}^P P, \end{cases} \tag{5}$$

where $\mathcal{D}^V$ and $\mathcal{D}^P$ arise from the discretization of $\partial v/\partial x$ and $\partial p/\partial x$, respectively. The time-dependent vectors $P$ and $V$ are the discrete approximations to the solution variables $p$ and $v$, respectively, and are of sizes $N^P$ and $N^V$, respectively. Matrices $\mathcal{A}^P$ and $\mathcal{A}^V$ are referred to as the norm matrices in the SBP literature. We limit ourselves to the case of diagonal norm matrices in this article.

The appearances of $\mathcal{A}^P$ and $\mathcal{A}^V$ in (5) may seem redundant at first sight, but are in fact important for the upcoming discussions. Loosely speaking, diagonal entries of $\mathcal{A}^P$ and $\mathcal{A}^V$ represent the areas that their corresponding grid points occupy. In [23], the authors pointed out that diagonal entries of the norm matrices and their corresponding discretization grid points provide quadrature rules for the underlying interval, acting as quadrature weights and quadrature points, respectively. This result will become useful in later sections where we derive the 2D SBP discretization operators.

To analyze the dynamical behavior of (5) with the energy method, we associate it with the discrete energy $E$, defined as

$$E = \frac{1}{2} P^T \mathcal{A}^P P + \frac{1}{2} V^T \mathcal{A}^V V. \tag{6}$$

Moreover, we define matrices $Q^V = \mathcal{A}^P \mathcal{D}^V$, $Q^P = \mathcal{A}^V \mathcal{D}^P$ and $Q = Q^V + \left(Q^P\right)^T$ to assist the upcoming discussion. Taking the temporal derivative of $E$ and substituting in (5), we obtain

$$\frac{dE}{dt} = P^T \mathcal{A}^P \frac{dP}{dt} + V^T \mathcal{A}^V \frac{dV}{dt} = -P^T \mathcal{A}^P \mathcal{D}^V V - V^T \mathcal{A}^V \mathcal{D}^P P = -P^T Q V. \tag{7}$$

We want to choose $Q$ properly such that $P^T Q V$ is the discrete correspondent of $(pv)\Big|_{x_L}^{x_R}$. In the existing literature, $Q$ is often chosen as

$$Q = \begin{bmatrix} -1 & & \\ & & \\ & & 1 \end{bmatrix}, \tag{8}$$

which leads to $-P^T Q V = P(1) \cdot V(1) - P(N^P) \cdot V(N^V)$. However, this quantity is not a good imitation to $(pv)\Big|_{x_L}^{x_R}$ for our discretization scenario since due to grid staggering, cf. Figure 4, $V(1)$ and $V(N^V)$ do not correspond to the evaluation of $v$ on the boundary grid points. To circumvent this issue, we ask $Q$ to have the following structure:

$$Q = -e_L \left(\mathcal{P}_L^V\right)^T + e_R \left(\mathcal{P}_R^V\right)^T, \tag{9}$$

where $e_L$ and $e_R$ are canonical basis vectors of the same size of $P$, taking the forms

$$e_L = [1, 0, \cdots, 0, 0]^T \quad \text{and} \quad e_R = [0, 0, \cdots, 0, 1]^T, \tag{10}$$

respectively. On the other hand, $\mathcal{P}_L^V$ and $\mathcal{P}_R^V$ are the discrete boundary projection operators, which are column vectors of the same size of $V$. When applied to $V$, $(\mathcal{P}_L^V)^T V$ and $(\mathcal{P}_R^V)^T V$ provide approximations to $v(x_L)$ and $v(x_R)$, respectively. With these definitions, we have

$$\frac{dE}{dt} \;=\; -P^T Q V \;=\; P(1) \cdot \left(\left(\mathcal{P}_L^V\right)^T V\right) - P(N^P) \cdot \left(\left(\mathcal{P}_R^V\right)^T V\right), \tag{11}$$

which is mimetic to $-(pv)\Big|_{x_L}^{x_R}$.

Several other discretization choices are presented in the following before we proceed to the resultant discretization matrices. First, we use the fourth-order standard staggered grid central difference stencil $[1/24,\ -9/8,\ 9/8,\ -1/24]/\Delta x$ to approximate $\partial/\partial x$ at the interior grid points. Second, we ask matrices $\mathcal{A}^P$ and $\mathcal{A}^V$ to be diagonal. According to [24], this implies that only stencils up to second-order accuracy can be derived for regions near the boundary. Third, the special regions near the boundary include four $p$ grid points and three $v$ grid points. We ask the derived stencil to be second-order accurate on these grid points. Finally, we ask second-order accuracy in the boundary projection operators $\mathcal{P}_L^V$ and $\mathcal{P}_R^V$.

We aim to design discretization matrices that satisfy the aforementioned accuracy constraints and lead to matrix $Q$ that has the structure presented in (9). Following the procedure summarized in [16], we arrive at the discretization matrices in (12a-12d) for problem (2) discretized on the grid shown in Figure 4. We note here that these matrices are presented for the unit grid spacing case, i.e., $\Delta x = 1$. When applied to general cases, $\mathcal{D}^V$ and $\mathcal{D}^P$ need to be scaled by dividing $\Delta x$; $\mathcal{A}^P$ and $\mathcal{A}^V$ need to be scaled by multiplying $\Delta x$.

$$
\mathcal{D}^V = \begin{bmatrix}
-2 & 3 & -1 & 0 & 0 & & & & & & & & \\
-1 & 1 & 0 & 0 & 0 & & & & & & & & \\
1/24 & -9/8 & 9/8 & -1/24 & 0 & & & & & & & & \\
-1/71 & 6/71 & -83/71 & 81/71 & -3/71 & & & & & & & & \\
 & & 1/24 & -9/8 & 9/8 & -1/24 & & & & & & & \\
 & & & 1/24 & -9/8 & 9/8 & -1/24 & & & & & & \\
 & & & & \ddots & \ddots & \ddots & \ddots & & & & & \\
 & & & & & & 1/24 & -9/8 & 9/8 & -1/24 & & & \\
 & & & & & & & 1/24 & -9/8 & 9/8 & -1/24 & & \\
 & & & & & & & & 3/71 & -81/71 & 83/71 & -6/71 & 1/71 \\
 & & & & & & & & 0 & 1/24 & -9/8 & 9/8 & -1/24 \\
 & & & & & & & & 0 & 0 & 0 & -1 & 1 \\
 & & & & & & & & 0 & 0 & 1 & -3 & 2
\end{bmatrix}; \tag{12a}
$$

$$
\mathcal{D}^P = \begin{bmatrix}
-79/78 & 27/26 & -1/26 & 1/78 & 0 & & & & & & & & \\
2/21 & -9/7 & 9/7 & -2/21 & 0 & & & & & & & & \\
1/75 & 0 & -27/25 & 83/75 & -1/25 & & & & & & & & \\
 & & 1/24 & -9/8 & 9/8 & -1/24 & & & & & & & \\
 & & & 1/24 & -9/8 & 9/8 & -1/24 & & & & & & \\
 & & & & \ddots & \ddots & \ddots & \ddots & & & & & \\
 & & & & & & 1/24 & -9/8 & 9/8 & -1/24 & & & \\
 & & & & & & & 1/24 & -9/8 & 9/8 & -1/24 & & \\
 & & & & & & & & 1/25 & -83/75 & 27/25 & 0 & -1/75 \\
 & & & & & & & & 0 & 2/21 & -9/7 & 9/7 & -2/21 \\
 & & & & & & & & 0 & -1/78 & 1/26 & -27/26 & 79/78
\end{bmatrix}; \tag{12b}
$$

$$
\mathcal{A}^P = \begin{bmatrix}
7/18 & & & & & & & & & & & & \\
 & 9/8 & & & & & & & & & & & \\
 & & 1 & & & & & & & & & & \\
 & & & 71/72 & & & & & & & & & \\
 & & & & 1 & & & & & & & & \\
 & & & & & 1 & & & & & & & \\
 & & & & & & \ddots & & & & & & \\
 & & & & & & & 1 & & & & & \\
 & & & & & & & & 1 & & & & \\
 & & & & & & & & & 71/72 & & & \\
 & & & & & & & & & & 1 & & \\
 & & & & & & & & & & & 9/8 & \\
 & & & & & & & & & & & & 7/18
\end{bmatrix}; \tag{12c}
$$

$$
\mathcal{A}^V = \begin{bmatrix}
13/12 & & & & & & & & & & \\
 & 7/8 & & & & & & & & & \\
 & & 25/24 & & & & & & & & \\
 & & & 1 & & & & & & & \\
 & & & & 1 & & & & & & \\
 & & & & & \ddots & & & & & \\
 & & & & & & 1 & & & & \\
 & & & & & & & 1 & & & \\
 & & & & & & & & 25/24 & & \\
 & & & & & & & & & 7/8 & \\
 & & & & & & & & & & 13/12
\end{bmatrix}. \tag{12d}
$$

These matrices lead to the following $Q$ matrix:

$$
Q = \begin{bmatrix}
-15/8 & 5/4 & -3/8 & & & \\
 & & & & & \\
 & & & & & \\
 & & & 3/8 & -5/4 & 15/8
\end{bmatrix}, \tag{13}
$$

while the boundary projection operators $\mathcal{P}_L^V$ and $\mathcal{P}_R^V$ are simply the transposes of the first and the last row of $Q$, respectively, which can be formally written as

$$\mathcal{P}_L^V = -\left(e_L^T Q\right)^T \quad \text{and} \quad \mathcal{P}_R^V = \left(e_R^T Q\right)^T. \tag{14}$$

**Remark 1.** *Here, we only address the case that p subgrid is aligned with both boundary points, cf. Figure 4, since this grid configuration is natural to impose free surface boundary condition, i.e.,* $p = 0$*, which is pertinent to the seismic applications in mind. However, if needed, other grid alignments at the boundary are possible and the corresponding discretization matrices can be derived in similar manner.*

*3.1.2. Imposing free surface boundary condition for the 1D model problem*

The technique of Simultaneous Approximation Terms (SATs), originally proposed in [15], can be used to impose various boundary conditions for various PDE systems in ways that lead to accurate and stable semi-discretized dynamical systems. Readers may refer to [25, 26] and the references therein for more information. In the following, we demonstrate how SATs can be used to weakly impose the free surface boundary condition for the 1D model problem described in (2). Specifically, PDE system (2) is closed with boundary conditions

$$p(x_L) = 0; \qquad p(x_R) = 0. \tag{15}$$

To impose these boundary conditions, we modify the semi-discretized system (5) by appending appropriate penalty terms, leading to the following system:

$$\begin{cases} \mathcal{A}^P \dfrac{dP}{dt} & = & -\mathcal{A}^P \mathcal{D}^V V; \\ \mathcal{A}^V \dfrac{dV}{dt} & = & -\mathcal{A}^V \mathcal{D}^P P + \sigma_L \cdot \mathcal{P}_L^V \cdot \left(e_L^T P - 0\right) + \sigma_R \cdot \mathcal{P}_R^V \cdot \left(e_R^T P - 0\right), \end{cases} \tag{16}$$

where $\sigma_L$ and $\sigma_R$ are scalar parameters to be chosen. The terms $\sigma_L \mathcal{P}_L^V \cdot (e_L^T P - 0)$ and $\sigma_R \mathcal{P}_L^R \cdot (e_R^T P - 0)$ of (16) are what the SATs refer to.

We now demonstrate that with proper choices of $\sigma_L$ and $\sigma_R$, the semi-discretized system (16) conserves the discrete energy given in (6). Going through the same energy analysis as in the previous section and substituting in (9), we obtain

$$\begin{aligned} \frac{dE}{dt} &= -P^T Q V + \sigma_L \left(V^T \mathcal{P}_L^V\right) \cdot \left(e_L^T P\right) + \sigma_R \left(V^T \mathcal{P}_L^R\right) \cdot \left(e_R^T P\right) \\ &= (\sigma_L + 1) \cdot \left(V^T \mathcal{P}_L^V\right) \cdot \left(e_L^T P\right) + (\sigma_R - 1) \cdot \left(V^T \mathcal{P}_R^V\right) \cdot \left(e_R^T P\right). \end{aligned} \tag{17}$$

We observe from (17) that by choosing $\sigma_L = -1$ and $\sigma_R = 1$, we have $dE/dt = 0$, i.e., the semi-discretized system (16) conserves energy $E$.

*3.1.3. Imposing interface conditions for the 1D model problem*

In this section, we consider a slightly different discretization scenario for the 1D model problem (2), where the 1D uniform grid is split into two segments by introducing an interface ($x_I$) to the $p$ subgrid, as demonstrated in Figure 5. The grid point at the interface is duplicated so that both segments own a grid point aligned at the interface. We attach superscripts $^-$ and $^+$ to solution variables $p$ and $v$, as well as their spatially discretized counterparts $P$ and $V$, in order to indicate to which side of the interface they belong. At the interface, we seek to impose the following interface conditions

$$p^- = p^+; \qquad v^- = v^+ \tag{18}$$

such that the wave can pass through the interface smoothly without interference. In this section, we explain how to impose the interface conditions in (18) with the SATs.

In order to focus on the interface treatment, we assume that both boundaries ($x_L$ and $x_R$) have been satisfactorily dealt with by some SBP discretization operators and SATs in an energy-conserving manner so that in the upcoming energy analysis, no term pops out in $Q$ due to the boundaries. From (12), we readily have the SBP discretization matrices for each individual segments. We use $\mathcal{A}^{P^-}$, $\mathcal{A}^{P^+}$, $\mathcal{A}^{V^-}$, $\mathcal{A}^{V^+}$, $\mathcal{D}^{V^-}$, $\mathcal{D}^{V^+}$, $\mathcal{D}^{P^-}$ and $\mathcal{D}^{P^+}$ to denote these discretization matrices, which have the same meaning as their counterparts in (5). Furthermore, we make the following definitions to simplify the upcoming discussion:

$$P = \begin{bmatrix} P^- \\ P^+ \end{bmatrix}, \; \mathcal{A}^P = \begin{bmatrix} \mathcal{A}^{P^-} & \mathbf{0} \\ \mathbf{0} & \mathcal{A}^{P^+} \end{bmatrix} \text{ and } \mathcal{D}^V = \begin{bmatrix} \mathcal{D}^{V^-} & \mathbf{0} \\ \mathbf{0} & \mathcal{D}^{V^+} \end{bmatrix}; \tag{19a}$$

$$V = \begin{bmatrix} V^- \\ V^+ \end{bmatrix}, \; \mathcal{A}^V = \begin{bmatrix} \mathcal{A}^{V^-} & \mathbf{0} \\ \mathbf{0} & \mathcal{A}^{V^+} \end{bmatrix} \text{ and } \mathcal{D}^P = \begin{bmatrix} \mathcal{D}^{P^-} & \mathbf{0} \\ \mathbf{0} & \mathcal{D}^{P^+} \end{bmatrix}. \tag{19b}$$

With these definitions, system (5) can still be used to describe the spatial discretization of PDE system (2) on the grid illustrated in Figure 5. We can now go through the same discrete energy analysis as in the previous section using the discrete energy defined in (6), and arrive again at relation

$$\frac{dE}{dt} = -P^T Q V, \tag{20}$$

except that $Q$ now takes the following form:

$$Q = \left[\begin{array}{ccc|ccc} & & & & & \\ & & & & & \\ {}^{3}\!/_{8} & -{}^{5}\!/_{4} & {}^{15}\!/_{8} & & & \\ \hline & & & -{}^{15}\!/_{8} & {}^{5}\!/_{4} & -{}^{3}\!/_{8} \\ & & & & & \\ & & & & & \end{array}\right]. \tag{21}$$

We note here that only two rows of $Q$ are nonzero, corresponding to the last entry of $P^-$ and the first entry of $P^+$.

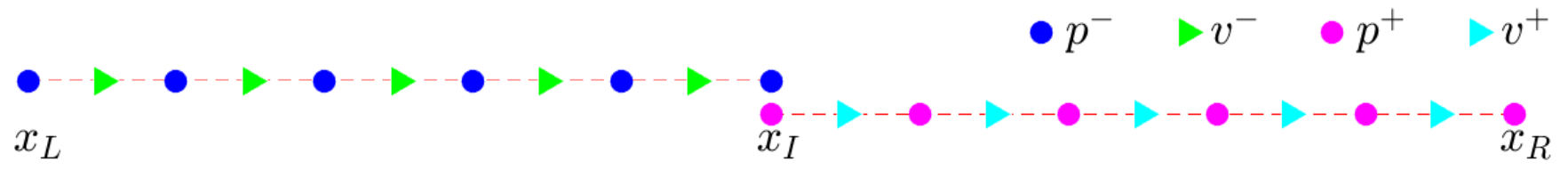

Figure 5: Illustration of the 1D staggered grid with an interface located on the $p$ subgrid.

Similar to $e_L$, $e_R$, $\mathcal{P}_L^V$ and $\mathcal{P}_R^V$ in Section 3.1.1, we introduce canonical basis vectors $e_I^-$ and $e_I^+$, and discrete interface projection operators $\mathcal{P}_I^{V^-}$ and $\mathcal{P}_I^{V^+}$ to simplify the discussion. Both $e_I^-$ and $e_I^+$ are of the same size of $P$, taking the forms

$$e_I^- = [0, \cdots, 1, 0, \cdots, 0]^T \quad \text{and} \quad e_I^+ = [0, \cdots, 0, 1, \cdots, 0]^T, \tag{22}$$

respectively. The only nonzero entry in $e_I^-$ corresponds to the last entry of $P^-$ while the only nonzero entry in $e_I^+$ corresponds to the first entry of $P^+$. Both $\mathcal{P}_I^{V^-}$ and $\mathcal{P}_I^{V^+}$ are column vectors of the same size of $V$, which can be formally defined as

$$\mathcal{P}_I^{V^-} = \left( (e_I^-)^T Q \right)^T \quad \text{and} \quad \mathcal{P}_I^{V^+} = -\left( (e_I^+)^T Q \right)^T, \tag{23}$$

respectively. When applied to $V$, $\left(\mathcal{P}_I^{V^-}\right)^T V$ and $\left(\mathcal{P}_I^{V^+}\right)^T V$ provide approximations to $v(x_I)$ from its left and from its right, respectively. With these notations, $Q$ in (21) can now be written succinctly as

$$Q = e_I^- \left(\mathcal{P}_I^{V^-}\right)^T - e_I^+ \left(\mathcal{P}_I^{V^+}\right)^T. \tag{24}$$

We proceed to show how the SATs technique can be used to impose the interface conditions in (18) in an energy-conserving manner. Specifically, we modify the semi-discretized system (5) as follows

$$
\begin{cases}
\mathcal{A}^P \dfrac{\partial P}{\partial t} & = & -\mathcal{A}^P \mathcal{D}^V V & + & \tau_I^- \cdot \mathbf{e}_I^- \cdot \left( \left( \mathcal{P}_I^{V^+} \right)^T V - \left( \mathcal{P}_I^{V^-} \right)^T V \right) \\
 & & & + & \tau_I^+ \cdot \mathbf{e}_I^+ \cdot \left( \left( \mathcal{P}_I^{V^+} \right)^T V - \left( \mathcal{P}_I^{V^-} \right)^T V \right); \\
\mathcal{A}^V \dfrac{\partial V}{\partial t} & = & -\mathcal{A}^V \mathcal{D}^P P & + & \sigma_I^- \cdot \left( \mathcal{P}_I^{V^-} \right) \cdot \left( (\mathbf{e}_I^+)^T P - (\mathbf{e}_I^-)^T P \right) \\
 & & & + & \sigma_I^+ \cdot \left( \mathcal{P}_I^{V^+} \right) \cdot \left( (\mathbf{e}_I^+)^T P - (\mathbf{e}_I^-)^T P \right),
\end{cases}
\tag{25}
$$

where $\tau_I^-$, $\tau_I^+$, $\sigma_I^-$ and $\sigma_I^+$ are scalar parameters to be chosen. Going through the same discrete energy analysis again for system (25), we arrive at

$$
\begin{aligned}
\frac{dE}{dt} &= (-1 - \tau_I^- - \sigma_I^-) \cdot (P^T e_I^-) \cdot \left( \left( \mathcal{P}_I^{V^-} \right)^T V \right) + (\tau_I^- - \sigma_I^+) \cdot (P^T e_I^-) \cdot \left( \left( \mathcal{P}_I^{V^+} \right)^T V \right) \\
&+ \left( -\tau_I^+ + \sigma_I^- \right) \cdot (P^T e_I^+) \cdot \left( \left( \mathcal{P}_I^{V^-} \right)^T V \right) + \left( 1 + \tau_I^+ + \sigma_I^+ \right) \cdot (P^T e_I^+) \cdot \left( \left( \mathcal{P}_I^{V^+} \right)^T V \right).
\end{aligned}
\tag{26}
$$

We observe that by choosing the set of parameters $\tau_I^-$, $\tau_I^+$, $\sigma_I^-$ and $\sigma_I^+$ properly so that linear system

$$
\begin{cases}
-1 - \tau_I^- - \sigma_I^- & = & 0 \\
\tau_I^- - \sigma_I^+ & = & 0 \\
-\tau_I^+ + \sigma_I^- & = & 0 \\
1 + \tau_I^+ + \sigma_I^+ & = & 0
\end{cases}
\tag{27}
$$

holds, we have ${dE}/{dt} = 0$ and consequently, system (25) conserves the discrete energy defined in (6). Linear system (27) is underdetermined and therefore has infinite solutions. For instance, $\tau_I^- = \tau_I^+ = \sigma_I^- = \sigma_I^+ = -1/2$ is one such solution.

### *3.2. The 2D acoustic wave equation*

In this section, we explain how the results of the previous section, including the discretization operators and the energy-conserving property, can be generalized to the 2D acoustic wave equation. To start, we consider the homogeneous version of (1) with unit constant coefficients:

$$
\begin{cases}
\dfrac{\partial p}{\partial t} & = & -\dfrac{\partial u}{\partial x} - \dfrac{\partial v}{\partial y}; \\
\dfrac{\partial u}{\partial t} & = & -\dfrac{\partial p}{\partial x}; \\
\dfrac{\partial v}{\partial t} & = & -\dfrac{\partial p}{\partial y},
\end{cases}
\tag{28}
$$

defined over rectangular domain $\square = (x_L, x_R) \times (y_L, y_R)$. Similar to the 1D case, we omit boundary conditions at this stage. Moreover, we define the energy associated with (28) as

$$
\mathscr{E} = \frac{1}{2} \int_{\square} \left( p^2 + u^2 + v^2 \right) d\square. \tag{29}
$$

Taking the temporal derivative of $\mathscr{E}$, substituting in (28) and then applying the divergence theorem, we arrive at

$$
\frac{d\mathscr{E}}{dt} = -\oint_{\partial\square} \left( p\vec{\mathbf{v}} \right) \cdot \vec{\mathbf{n}}\, d\partial_{\square} = \int_{\partial B} pv\, dx - \int_{\partial T} pv\, dx + \int_{\partial L} pu\, dy - \int_{\partial R} pu\, dy, \tag{30}
$$

where vector field $\vec{\mathbf{v}}$ stands for $[u, v]^T$, $\vec{\mathbf{n}}$ stands for the outward-pointing unit normal vector, $\partial_{\square}$ stands for the entire boundary of domain $\square$, and $\partial B$, $\partial T$, $\partial L$ and $\partial R$ stand for the bottom, top, left and right boundaries of $\square$, respectively.

Similarly to the 1D case, we observe here that the energy evolution depends only on the solution at the boundary. We aim to retain a similar property in the semi-discretized system.

In the following, we extend the notational system to ease the upcoming discussion. To start, we consider a uniform staggered grid as illustrated in Figure 1. Each 2D subgrid is built as the tensor product of its corresponding 1D subgrids. Taking variable $p$ as an example, we use $\mathbf{x}^P = \left[x_1^P, \cdots, x_{N_x^P}^P\right]$ and $\mathbf{y}^P = \left[y_1^P, \cdots, y_{N_y^P}^P\right]$ to denote the 1D subgrids that discretizes intervals $[x_L, x_R]$ and $[y_L, y_R]$, respectively, where $x_1^P = x_L$, $x_{N_x^P}^P = x_R$, $y_1^P = y_L$ and $y_{N_y^P}^P = y_R$. Similarly, we associate variable $u$ with 1D subgrids $\mathbf{x}^U$ and $\mathbf{y}^U$, and associate variable $v$ with 1D subgrids $\mathbf{x}^V$ and $\mathbf{y}^V$. Due to grid staggering on the respective directions, $\mathbf{x}^U$ and $\mathbf{y}^V$ consist of midpoints of those grid points in $\mathbf{x}^P$ and $\mathbf{y}^P$, respectively. Meanwhile, $\mathbf{x}^V$ and $\mathbf{y}^U$ are identical to $\mathbf{x}^P$ and $\mathbf{y}^P$, respectively.

Moreover, we use column vectors $P$, $U$ and $V$ to denote the spatially discretized solution variables. They are mapped from the respective 2D subgrids in column-wise manner. Taking variable $p$ as an example, we use $P_i^y$ to denote the restriction of $P$ on the $i$th grid column for $i = 1, \cdots, N_x^P$.

Next, we introduce notations for 1D SBP discretization matrices of two relevant PDE systems. These 1D matrices are building blocks for the 2D SBP discretization matrices presented in Section 3.2.1. Specifically, we use $\underline{\mathbf{a}}_y$, $\overline{\mathbf{a}}_y$, $\mathbf{d}_y^V$ and $\mathbf{d}_y^P$ to denote the 1D SBP discretization matrices for PDE system

$$\begin{cases} \dfrac{\partial p^y}{\partial t} &= -\dfrac{\partial v^y}{\partial y}; \\ \dfrac{\partial v^y}{\partial t} &= -\dfrac{\partial p^y}{\partial y}, \end{cases} \tag{31}$$

defined over $(y_L, y_R)$ and discretized on the staggered grid composed of $\mathbf{y}^P$ and $\mathbf{y}^V$. Matrices $\underline{\mathbf{a}}_y$, $\overline{\mathbf{a}}_y$, $\mathbf{d}_y^V$ and $\mathbf{d}_y^P$ have similar meanings to their counterparts $\mathcal{A}^P$, $\mathcal{A}^V$, $\mathcal{D}^V$ and $\mathcal{D}^P$ in Section 3.1.1, in that order. We associate symbol $\underline{\mathbf{a}}_y$ (with a bar underlying $\mathbf{a}$) with $\mathbf{y}^P$ and symbol $\overline{\mathbf{a}}_y$ (with a bar overlying $\mathbf{a}$) with $\mathbf{y}^V$. Following this naming convention, we use symbols $\underline{I}_y$ and $\overline{I}_y$ to denote the identity matrices having the same sizes of $\underline{\mathbf{a}}_y$ and $\overline{\mathbf{a}}_y$, respectively. Furthermore, we define matrices $\mathbf{q}_y^V$, $\mathbf{q}_y^P$ and $\mathbf{q}_y$ as

$$\mathbf{q}_y^V = \underline{\mathbf{a}}_y \mathbf{d}_y^V, \quad \mathbf{q}_y^P = \overline{\mathbf{a}}_y \mathbf{d}_y^P \quad \text{and} \quad \mathbf{q}_y = \mathbf{q}_y^V + \left(\mathbf{q}_y^P\right)^T. \tag{32}$$

Given the SBP property of the discretization matrices, $\mathbf{q}_y$ can be written as

$$\mathbf{q}_y = -e_{y_L}\left(\mathcal{P}_{y_L}^V\right)^T + e_{y_R}\left(\mathcal{P}_{y_R}^V\right)^T, \tag{33}$$

where $e_{y_L}$, $e_{y_R}$, $\mathcal{P}_{y_L}^V$ and $\mathcal{P}_{y_R}^V$ have similar meanings to their counterparts in Section 3.1.1.

On the other hand, we use $|\mathbf{a}_x$, $\mathbf{a}|_x$, $\mathbf{d}_x^U$ and $\mathbf{d}_x^P$ to denote the 1D SBP discretization matrices for PDE system

$$\begin{cases} \dfrac{\partial p^x}{\partial t} &= -\dfrac{\partial u^x}{\partial x}; \\ \dfrac{\partial u^x}{\partial t} &= -\dfrac{\partial p^x}{\partial x}, \end{cases} \tag{34}$$

defined over $(x_L, x_R)$ and discretized on the staggered grid composed of $\mathbf{x}^P$ and $\mathbf{x}^U$. Matrices $|\mathbf{a}_x$, $\mathbf{a}|_x$, $\mathbf{d}_x^U$ and $\mathbf{d}_x^P$ have similar meanings to their counterparts $\mathcal{A}^P$, $\mathcal{A}^V$, $\mathcal{D}^V$ and $\mathcal{D}^P$ in Section 3.1.1, in that order. We associate symbol $|\mathbf{a}_x$ (with a bar on the left of $\mathbf{a}$) with $\mathbf{x}^P$ and $\mathbf{a}|_x$ (with a bar on the right of $\mathbf{a}$) with $\mathbf{x}^U$. Following this naming convention, we use symbols $|I_y$ and $I|_y$ to denote the identity matrices having the same sizes of $|\mathbf{a}_x$ and $\mathbf{a}|_x$, respectively. Furthermore, we define matrices $\mathbf{q}_x^U$, $\mathbf{q}_x^P$ and $\mathbf{q}_x$ as

$$\mathbf{q}_x^U = |\mathbf{a}_x \mathbf{d}_x^U, \quad \mathbf{q}_x^P = \mathbf{a}|_x \mathbf{d}_x^P \quad \text{and} \quad \mathbf{q}_x = \mathbf{q}_x^U + \left(\mathbf{q}_x^P\right)^T. \tag{35}$$

Given the SBP property of the discretization matrices, $\mathbf{q}_x$ can be written as

$$\mathbf{q}_x = -e_{x_L}\left(\mathcal{P}_{x_L}^U\right)^T + e_{x_R}\left(\mathcal{P}_{x_R}^U\right)^T, \tag{36}$$

where $e_{x_L}$, $e_{x_R}$, $\mathcal{P}_{x_L}^U$ and $\mathcal{P}_{x_R}^U$ have similar meanings to their counterparts in Section 3.1.1.

### 3.2.1. SBP discretization operators built via tensor product

Using these 1D SBP discretization matrices as building blocks, we construct the 2D discretization matrices for PDE system (28) as:

$$
\begin{aligned}
&\mathcal{A}^P = |\mathbf{a}_x \otimes \underline{\mathbf{a}}_y, \quad \mathcal{A}^U = \mathbf{a}|_x \otimes \underline{\mathbf{a}}_y, \quad \mathcal{A}^V = |\mathbf{a}_x \otimes \overline{\mathbf{a}}_y, \\
&\mathcal{D}_y^V = |I_x \otimes \mathbf{d}_y^V, \quad \mathcal{D}_y^P = |I_x \otimes \mathbf{d}_y^P, \quad \mathcal{D}_x^U = \mathbf{d}_x^U \otimes \underline{I}_y, \quad \mathcal{D}_x^P = \mathbf{d}_x^P \otimes \underline{I}_y,
\end{aligned}
\tag{37}
$$

where the symbol '⊗' stands for the tensor product operation on matrices. Readers may refer to [27, 28] for the definition, properties and applications of tensor product of matrices.

Matrices $\mathcal{D}_y^V$, $\mathcal{D}_y^P$, $\mathcal{D}_x^U$ and $\mathcal{D}_x^P$ can be applied on $V$, $P$, $U$, $P$, respectively, providing approximations of $dv/dy$, $dp/dy$, $du/dx$ and $dp/dx$ to orders of accuracy consistent with those of $\mathbf{d}_y^V$, $\mathbf{d}_y^P$, $\mathbf{d}_x^U$ and $\mathbf{d}_x^P$, respectively. Convergence tests concerning accuracy of these operators can be found in Appendix A, where both discretization scenarios of Section 3.2.2 (without interface) and Section 3.2.3 (with interface) are considered.

These matrices lead to the following semi-discretization of PDE system (28):

$$
\begin{cases}
\mathcal{A}^P \dfrac{dP}{dt} &= -\mathcal{A}^P \mathcal{D}_x^U U - \mathcal{A}^P \mathcal{D}_y^V V; \\
\mathcal{A}^U \dfrac{dU}{dt} &= -\mathcal{A}^U \mathcal{D}_x^P P; \\
\mathcal{A}^V \dfrac{dV}{dt} &= -\mathcal{A}^V \mathcal{D}_y^P P.
\end{cases}
\tag{38}
$$

We define the discrete energy associated with (38) as

$$
E = \frac{1}{2} P^T \mathcal{A}^P P + \frac{1}{2} U^T \mathcal{A}^U U + \frac{1}{2} V^T \mathcal{A}^V V.
\tag{39}
$$

Furthermore, we define matrices $Q_y^V$, $Q_y^P$, $Q_x^U$ and $Q_x^P$ as

$$
Q_y^V = \mathcal{A}^P \mathcal{D}_y^V, \quad Q_y^P = \mathcal{A}^V \mathcal{D}_y^P, \quad Q_x^U = \mathcal{A}^P \mathcal{D}_x^U \quad \text{and} \quad Q_x^P = \mathcal{A}^U \mathcal{D}_x^P,
\tag{40}
$$

respectively. Substituting the definitions in (37), (32) and (35) into (40) and applying the mixed-product rule of tensor product, we have

$$
Q_y^V = |\mathbf{a}_x \otimes \mathbf{q}_y^V, \quad Q_y^P = |\mathbf{a}_x \otimes \mathbf{q}_y^P, \quad Q_x^U = \mathbf{q}_x^U \otimes \underline{\mathbf{a}}_y, \quad Q_x^P = \mathbf{q}_x^P \otimes \underline{\mathbf{a}}_y.
\tag{41}
$$

We proceed to define matrices $Q_y$ and $Q_x$ as

$$
Q_y = Q_y^V + \left(Q_y^P\right)^T \quad \text{and} \quad Q_x = Q_x^U + \left(Q_x^P\right)^T,
\tag{42}
$$

respectively. Recalling the definitions of $\mathbf{q}_y$ in (32) and $\mathbf{q}_x$ in (35) as well as relations (33) and (36), matrices $Q_y$ and $Q_x$ can be written as:

$$
Q_y = |\mathbf{a}_x \otimes \mathbf{q}_y = -|\mathbf{a}_x \otimes \left(e_{y_L} \left(\mathcal{P}_{y_L}^V\right)^T\right) + |\mathbf{a}_x \otimes \left(e_{y_R} \left(\mathcal{P}_{y_R}^V\right)^T\right);
\tag{43a}
$$

$$
Q_x = \mathbf{q}_x \otimes \underline{\mathbf{a}}_y = -\left(e_{x_L} \left(\mathcal{P}_{x_L}^U\right)^T\right) \otimes \underline{\mathbf{a}}_y + \left(e_{x_R} \left(\mathcal{P}_{x_R}^U\right)^T\right) \otimes \underline{\mathbf{a}}_y.
\tag{43b}
$$

With these notations and relations, we now proceed to the discrete energy analysis by taking the temporal derivative of (39) and then substituting in (38), obtaining

$$
\begin{aligned}
\frac{dE}{dt} &= -P^T Q_y V - P^T Q_x U \\
&= P^T \left(|\mathbf{a}_x \otimes \left(e_{y_L} \left(\mathcal{P}_{y_L}^V\right)^T\right)\right) V - P^T \left(|\mathbf{a}_x \otimes \left(e_{y_R} \left(\mathcal{P}_{y_R}^V\right)^T\right)\right) V \\
&+ P^T \left(\left(e_{x_L} \left(\mathcal{P}_{x_L}^U\right)^T\right) \otimes \underline{\mathbf{a}}_y\right) U - P^T \left(\left(e_{x_R} \left(\mathcal{P}_{x_R}^U\right)^T\right) \otimes \underline{\mathbf{a}}_y\right) U.
\end{aligned}
\tag{44}
$$

It can be shown that the four terms in the last relation of (44) are discrete correspondents of the four terms in the last relation of (30), respectively, in the listed order. This means that the discrete energy change rate $dE/dt$ is mimetic to its continuous counterpart $d\mathscr{E}/dt$, cf. (30), and therefore, semi-discretized system (38) is an SBP discretization of PDE system (28). In the interest of conserving space, we only demonstrate for the first term in the following.

Using the mixed-product rule of tensor product, we can write

$$|\mathbf{a}_x \otimes \left(e_{y_L}\left(\mathcal{P}^V_{y_L}\right)^T\right) = \left(|I_x \cdot |\mathbf{a}_x \cdot |I_x\right) \otimes \left(e_{y_L} \cdot 1 \cdot \left(\mathcal{P}^V_{y_L}\right)^T\right) = \left(|I_x \otimes e_{y_L}\right) \cdot |\mathbf{a}_x \cdot \left(|I_x \otimes \left(\mathcal{P}^V_{y_L}\right)^T\right),$$

where '1' stands for the 1 by 1 matrix with entry 1. Consequently, we have

$$P^T\left(|\mathbf{a}_x \otimes \left(e_{y_L}\left(\mathcal{P}^V_{y_L}\right)^T\right)\right)V = \left(P^T\left(|I_x \otimes e_{y_L}\right)\right) \cdot |\mathbf{a}_x \cdot \left(\left(|I_x \otimes \left(\mathcal{P}^V_{y_L}\right)^T\right)V\right). \tag{45}$$

Using the definition of tensor product, we can write the following two expressions:

$$P^T\left(|I_x \otimes e_{y_L}\right) = \left[\left(P^y_1\right)^T e_{y_L}, \cdots, \left(P^y_{N^P_x}\right)^T e_{y_L}\right]; \tag{46}$$

$$\left(|I_x \otimes \left(\mathcal{P}^V_{y_L}\right)^T\right)V = \left[\left(\mathcal{P}^V_{y_L}\right)^T V^y_1 \cdots \left(\mathcal{P}^V_{y_L}\right)^T V^y_{N^V_x}\right]^T. \tag{47}$$

Since each $P^y_i$ corresponds to a $p$ subgrid column, $\left(P^y_i\right)^T e_{y_L}$ approximates $p$ at grid point $\left(x^P_i, y^P_1\right)$, which is on the bottom boundary $\partial B$. Similarly, $\left(\mathcal{P}^V_{y_L}\right)^T V^y_i$ approximates $v$ at grid point $(x^V_i, y^P_1)$ on $\partial B$. Recalling that $v$ subgrid is not staggered in $x$-direction with respect to $p$ subgrid, we have $x^V_i = x^P_i$. Substituting (46) and (47) into (45) and realizing that $|\mathbf{a}_x$ is a diagonal matrix, we have

$$\begin{aligned} P^T\left(|\mathbf{a}_x \otimes \left(e_{y_L}\left(\mathcal{P}^V_{y_L}\right)^T\right)\right)V &= \sum_{i=1}^{N^P_x} |\mathbf{a}_x(i,i) \cdot \left(\left(P^y_i\right)^T e_{y_L}\right) \cdot \left(\left(\mathcal{P}^V_{y_L}\right)^T V^y_i\right) \\ &\approx \sum_{i=1}^{N^P_x} |\mathbf{a}_x(i,i) \cdot p(x^P_i, y^P_1) \cdot v(x^P_i, y^P_1). \end{aligned} \tag{48}$$

Recalling that the diagonal entries of $|\mathbf{a}_x$ define a set of quadrature weights for interval $[x_L, x_R]$, with the corresponding quadrature points given by $\mathbf{x}^P = \left[x^P_1, \cdots, x^P_{N^P_x}\right]$, we conclude that the term $P^T\left(|\mathbf{a}_x \otimes \left(e_{y_L}\left(\mathcal{P}^V_{y_L}\right)^T\right)\right)V$ in (44) is a discrete correspondence of the term $\int_{\partial B} pv\,dx$ in (30).

*3.2.2. Imposing free surface boundary condition for the 2D acoustic wave equation*

In this section, we briefly demonstrate how the free surface boundary condition $p = 0$ on $\partial\Box$ can be imposed via the SATs technique for PDE system (28). Given the notations and results in the previous section, the technique presented in Section 3.1.2 can be extended to this case naturally. Specifically, we append penalty terms to semi-discretized system (38), leading to:

$$\left\{\begin{aligned} \mathcal{A}^P\frac{dP}{dt} &= -\mathcal{A}^P\mathcal{D}^U_x U - \mathcal{A}^P\mathcal{D}^V_y V; \\ \mathcal{A}^U\frac{dU}{dt} &= -\mathcal{A}^U\mathcal{D}^P_x P + \sigma_L\left(\mathcal{P}^U_{x_L} \otimes \underline{I}_y\right) \cdot \underline{\mathbf{a}}_y \cdot \left(\left(\left(e_{x_L}\right)^T \otimes \underline{I}_y\right)P - \underline{\mathbf{0}}_y\right) \\ &\quad + \sigma_R\left(\mathcal{P}^U_{x_R} \otimes \underline{I}_y\right) \cdot \underline{\mathbf{a}}_y \cdot \left(\left(\left(e_{x_R}\right)^T \otimes \underline{I}_y\right)P - \underline{\mathbf{0}}_y\right); \\ \mathcal{A}^V\frac{dV}{dt} &= -\mathcal{A}^V\mathcal{D}^P_y P + \sigma_B\left(|I_x \otimes \mathcal{P}^V_{y_L}\right) \cdot |\mathbf{a}_x \cdot \left(\left(|I_x \otimes \left(e_{y_L}\right)^T\right)P - |\mathbf{0}_x\right) \\ &\quad + \sigma_T\left(|I_x \otimes \mathcal{P}^V_{y_R}\right) \cdot |\mathbf{a}_x \cdot \left(\left(|I_x \otimes \left(e_{y_R}\right)^T\right)P - |\mathbf{0}_x\right), \end{aligned}\right. \tag{49}$$

where $\sigma_L$, $\sigma_R$, $\sigma_B$ and $\sigma_T$ are parameters to be chosen, while $|\mathbf{0}_x$ and $\underline{\mathbf{0}}_y$ are zero vectors having the sizes of $\mathbf{x}^P$ and $\mathbf{y}^P$, respectively. Going through the energy analysis with the discrete energy given in (39), we find that by choosing $\sigma_L = -1$, $\sigma_R = 1$, $\sigma_B = -1$ and $\sigma_T = 1$, we have ${}^{dE}/_{dt} = 0$, i.e., system (49) conserves energy (38).

### *3.2.3. SATs for the 2D interface*

In this section, we consider a staggered block-wise uniform grid as illustrated in Figure 6. The entire simulation domain is split into two blocks by a horizontal interface aligned with the $p$ subgrid. The $p$ subgrid of both the top block and the bottom block have grid points (with different grid spacings) allocated on the interface. In order to distinguish, we attach superscript $^-$ to symbols associated with the bottom block and $^+$ to symbols associated with the top block.

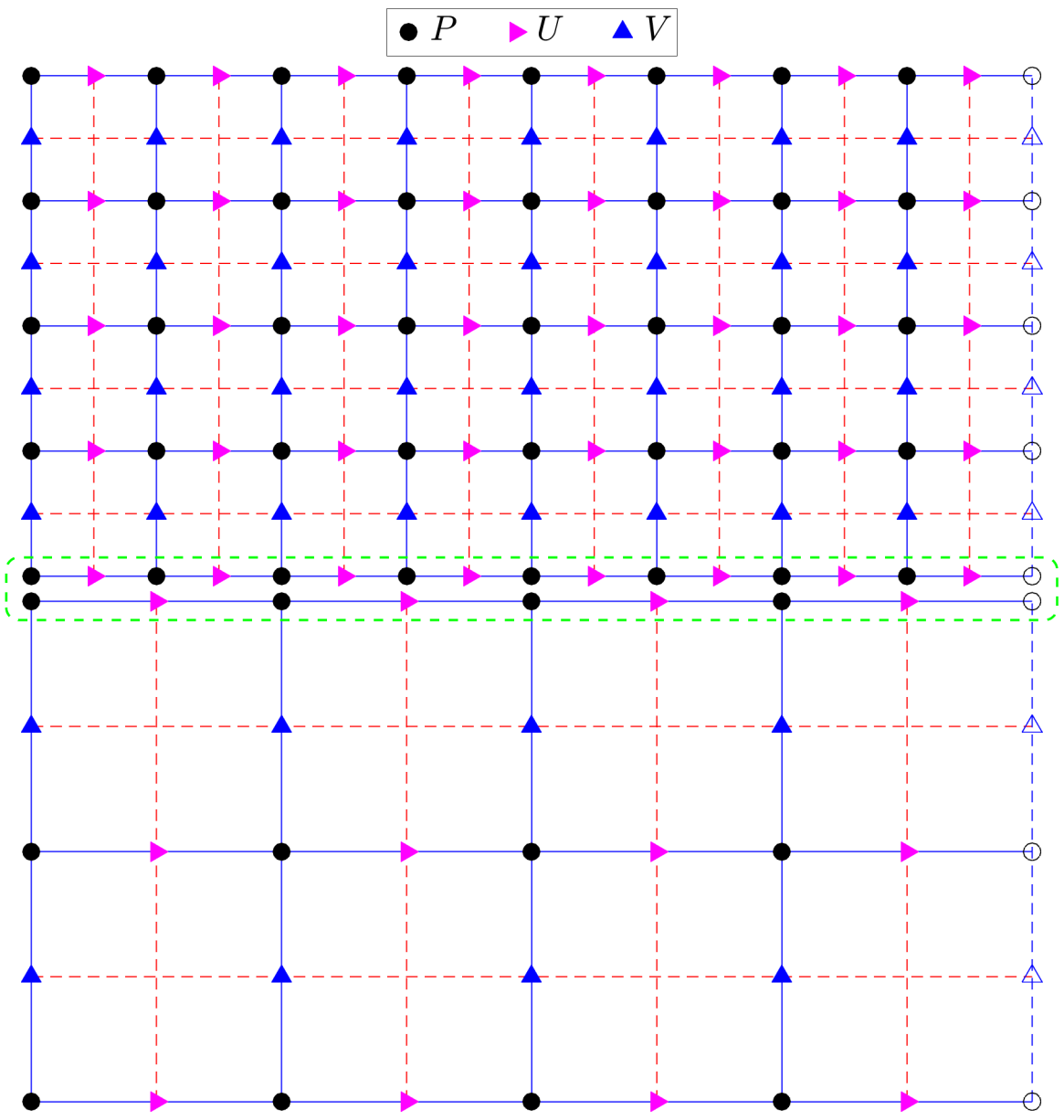

Figure 6: Illustration of a block-wise uniform grid consisting of two horizontal layers with a 2:1 ratio of the grid spacing. Each layer owns a grid row aligned along the interface.

Neglecting the boundary and interface conditions for now, on each block, PDE system (28) can be discretized with the SBP operators presented in Section 3.2.1, leading to the following semi-discretized systems:

$$\left\{\begin{aligned} \mathcal{A}^{P^-}\frac{dP^-}{dt} &= -\mathcal{A}^{P^-}\mathcal{D}_x^{U^-}U^- - \mathcal{A}^{P^-}\mathcal{D}_y^{V^-}V^-; \\ \mathcal{A}^{U^-}\frac{dU^-}{dt} &= -\mathcal{A}^{U^-}\mathcal{D}_x^{P^-}P^-; \\ \mathcal{A}^{V^-}\frac{dV^-}{dt} &= -\mathcal{A}^{V^-}\mathcal{D}_y^{P^-}P^-; \end{aligned}\right. \tag{50a}$$

$$\left\{\begin{aligned} \mathcal{A}^{P^+}\frac{dP^+}{dt} &= -\mathcal{A}^{P^+}\mathcal{D}_x^{U^+}U^+ - \mathcal{A}^{P^+}\mathcal{D}_y^{V^+}V^+; \\ \mathcal{A}^{U^+}\frac{dU^+}{dt} &= -\mathcal{A}^{U^+}\mathcal{D}_x^{P^+}P^+; \\ \mathcal{A}^{V^+}\frac{dV^+}{dt} &= -\mathcal{A}^{V^+}\mathcal{D}_y^{P^+}P^+. \end{aligned}\right. \tag{50b}$$

Symbols in (50a) and (50b) have the same meanings as their counterparts (without the superscripts) in Section 3.2.1.

The discrete energies associated with (50a) and (50b) are defined as

$$E^- = \frac{1}{2}(P^-)^T \mathcal{A}^{P^-} P^- + \frac{1}{2}(U^-)^T \mathcal{A}^{U^-} U^- + \frac{1}{2}(V^-)^T \mathcal{A}^{V^-} V^-; \tag{51a}$$

$$E^+ = \frac{1}{2}(P^+)^T \mathcal{A}^{P^+} P^+ + \frac{1}{2}(U^+)^T \mathcal{A}^{U^+} U^+ + \frac{1}{2}(V^+)^T \mathcal{A}^{V^+} V^+, \tag{51b}$$

respectively. Additionally, we define the total energy $E$ as the sum of $E^-$ and $E^+$.

Similarly to the 1D case presented in Section 3.1.3, we assume that all boundaries have been satisfactorily dealt with by some SBP operators and SATs in an energy-conserving manner so that in the upcoming discussion concerning the discrete energy analysis of (50), we only need to address terms related to the interface. Going through the same procedure as in Section 3.2.1 for systems (50a) and (50b), we arrive at:

$$\frac{dE^-}{dt} = -(P^-)^T \cdot \left(|I_x^- \otimes e_y^-\right) \cdot |\mathbf{a}_x^- \cdot \left(|I_x^- \otimes \left(\mathcal{P}_y^{V^-}\right)^T\right) V^-; \tag{52a}$$

$$\frac{dE^+}{dt} = \ (P^+)^T \cdot \left(|I_x^+ \otimes e_y^+\right) \cdot |\mathbf{a}_x^+ \cdot \left(|I_x^+ \otimes \left(\mathcal{P}_y^{V^+}\right)^T\right) V^+, \tag{52b}$$

respectively, where 1D column vectors $e_y^-$, $e_y^+$, $\mathcal{P}_y^{V^-}$ and $\mathcal{P}_y^{V^+}$ have similar meanings to their counterparts in Section 3.1.3. Specifically, $(e_y^-)^T P_i^{y^-}$ and $(e_y^+)^T P_i^{y^+}$ select values of the respective entries in $P_i^{y^-}$ and $P_i^{y^+}$ that correspond to the interface, respectively, while $(\mathcal{P}_y^{V^-})^T V_i^{y^-}$ and $(\mathcal{P}_y^{V^+})^T V_i^{y^+}$ project values in $V_i^{y^-}$ and $V_i^{y^+}$ to the interface, respectively.

We aim to couple the two semi-discretized systems (50a) and (50b) by properly imposing the following interface conditions:

$$p^- = p^+; \tag{53a}$$

$$\vec{\mathbf{v}}^- \cdot \vec{\mathbf{n}}^- + \vec{\mathbf{v}}^+ \cdot \vec{\mathbf{n}}^+ = 0, \tag{53b}$$

where vector fields $\vec{\mathbf{v}}^-$ and $\vec{\mathbf{v}}^+$ stand for $[u^-, v^-]^T$ and $[u^+, v^+]^T$, respectively, while $\vec{\mathbf{n}}^-$ and $\vec{\mathbf{n}}^+$ are the outward-pointing unit normal vectors on the interface. For the particular case under consideration where the interface is placed horizontally (cf. Figure 6), we have $\vec{\mathbf{n}}^- = [0, 1]^T$ and $\vec{\mathbf{n}}^+ = [0, -1]^T$, and therefore, (53b) boils down to

$$v^- = v^+. \tag{54}$$

We append systems (50a) and (50b) with penalty terms to account for the interface conditions in (53), leading to systems (55a) and (55b):

$$\begin{cases} \mathcal{A}^{P^-} \dfrac{dP^-}{dt} = -\mathcal{A}^{P^-} \mathcal{D}_x^{U^-} U^- - \mathcal{A}^{P^-} \mathcal{D}_y^{V^-} V^- \\ \qquad + \sigma_P^- \cdot \left(|I_x^- \otimes e_y^-\right) \cdot |\mathbf{a}_x^- \cdot \left[\mathcal{T}_+^- \left(\left(|I_x^+ \otimes \left(\mathcal{P}_y^{V^+}\right)^T\right) V^+\right) - \left(|I_x^- \otimes \left(\mathcal{P}_y^{V^-}\right)^T\right) V^-\right]; \\ \mathcal{A}^{U^-} \dfrac{dU^-}{dt} = -\mathcal{A}^{U^-} \mathcal{D}_x^{P^-} P^-; \\ \mathcal{A}^{V^-} \dfrac{dV^-}{dt} = -\mathcal{A}^{V^-} \mathcal{D}_y^{P^-} P^- \\ \qquad + \sigma_V^- \cdot \left(|I_x^- \otimes \mathcal{P}_y^{V^-}\right) \cdot |\mathbf{a}_x^- \cdot \left[\mathcal{T}_+^- \left(\left(|I_x^+ \otimes \left(e_y^+\right)^T\right) P^+\right) - \left(|I_x^- \otimes \left(e_y^-\right)^T\right) P^-\right]; \end{cases} \tag{55a}$$

$$\begin{cases} \mathcal{A}^{P^+} \dfrac{dP^+}{dt} = -\mathcal{A}^{P^+} \mathcal{D}_x^{U^+} U^+ - \mathcal{A}^{P^+} \mathcal{D}_y^{V^+} V^+ \\ \qquad + \sigma_P^+ \cdot \left(|I_x^+ \otimes e_y^+\right) \cdot |\mathbf{a}_x^+ \cdot \left[\left(|I_x^+ \otimes \left(\mathcal{P}_y^{V^+}\right)^T\right) V^+ - \mathcal{T}_-^+ \left(\left(|I_x^- \otimes \left(\mathcal{P}_y^{V^-}\right)^T\right) V^-\right)\right]; \\ \mathcal{A}^{U^+} \dfrac{dU^+}{dt} = -\mathcal{A}^{U^+} \mathcal{D}_x^{P^+} P^+; \\ \mathcal{A}^{V^+} \dfrac{dV^+}{dt} = -\mathcal{A}^{V^+} \mathcal{D}_y^{P^+} P^+ \\ \qquad + \sigma_V^+ \cdot \left(|I_x^+ \otimes \mathcal{P}_y^{V^+}\right) \cdot |\mathbf{a}_x^+ \cdot \left[\left(|I_x^+ \otimes \left(e_y^+\right)^T\right) P^+ - \mathcal{T}_-^+ \left(\left(|I_x^- \otimes \left(e_y^-\right)^T\right) P^-\right)\right], \end{cases} \tag{55b}$$

respectively, where $\sigma_P^-$, $\sigma_V^-$, $\sigma_P^+$ and $\sigma_V^+$ are scalar parameters to be chosen. To simplify the upcoming discussion, we define 1D column vectors $\mathbf{p}_{x_I}^-$, $\mathbf{p}_{x_I}^+$, $\mathbf{v}_{x_I}^-$ and $\mathbf{v}_{x_I}^+$ as:

$$\begin{aligned} \mathbf{p}_{x_I}^- &= \left( |I_x^- \otimes \left(e_y^-\right)^T \right) P^-, \quad & \mathbf{p}_{x_I}^+ &= \left( |I_x^+ \otimes \left(e_y^+\right)^T \right) P^+, \\ \mathbf{v}_{x_I}^- &= \left( |I_x^- \otimes \left(\mathcal{P}_y^{V^-}\right)^T \right) V^-, \quad & \mathbf{v}_{x_I}^+ &= \left( |I_x^+ \otimes \left(\mathcal{P}_y^{V^+}\right)^T \right) V^+. \end{aligned} \tag{56}$$

These 1D vectors can be understood as the restriction or projection of the respective 2D fields on the interface. In general, vectors $\mathbf{p}_{x_I}^-$ and $\mathbf{p}_{x_I}^+$ cannot be summed or subtracted from each other since the grids on the two sides of the interface are allowed to be nonconforming. The same holds for vectors $\mathbf{v}_{x_I}^-$ and $\mathbf{v}_{x_I}^+$. The interface transfer operators $\mathcal{T}_+^-$ and $\mathcal{T}_-^+$ are thus introduced. Specifically, $\mathcal{T}_+^-$ is an interpolation operator that applies on $\mathbf{p}_{x_I}^+$ and $\mathbf{v}_{x_I}^+$ so that the resulting vectors are compatible with $\mathbf{p}_{x_I}^-$ and $\mathbf{v}_{x_I}^-$ for sum or subtraction. Similarly, $\mathcal{T}_-^+$ is an interpolation operator that applies on $\mathbf{p}_{x_I}^-$ and $\mathbf{v}_{x_I}^-$ so that the resulting vectors are compatible with $\mathbf{p}_{x_I}^+$ and $\mathbf{v}_{x_I}^+$ for sum or subtraction. By design, we ask $\mathcal{T}_+^-$ and $\mathcal{T}_-^+$ to have the same order of accuracy as the projection operators $\mathcal{P}_y^{V^-}$ and $\mathcal{P}_y^{V^+}$.

Going through the discrete energy analysis for (55), we arrive at

$$\begin{aligned} \frac{dE}{dt} \;=\; &- \; (1+\sigma_P^- + \sigma_V^-)\cdot \mathbf{p}_{x_I}^- \cdot |\mathbf{a}_{x_I}^- \cdot \mathbf{v}_{x_I}^- \\ &+ \; (1+\sigma_P^+ + \sigma_V^+)\cdot \mathbf{p}_{x_I}^+ \cdot |\mathbf{a}_{x_I}^+ \cdot \mathbf{v}_{x_I}^+ \\ &- \; \mathbf{p}_{x_I}^- \cdot \left( \sigma_V^+ \left(\mathcal{T}_-^+\right)^T |\mathbf{a}_{x_I}^+ - \sigma_P^- |\mathbf{a}_{x_I}^- \mathcal{T}_+^- \right) \cdot \mathbf{v}_{x_I}^+ \\ &- \; \mathbf{p}_{x_I}^+ \cdot \left( \sigma_P^+ |\mathbf{a}_{x_I}^+ \mathcal{T}_-^+ - \sigma_V^- \left(\mathcal{T}_+^-\right)^T |\mathbf{a}_{x_I}^- \right) \cdot \mathbf{v}_{x_I}^-. \end{aligned} \tag{57}$$

Aside from the accuracy requirement, we impose the additional constraint:

$$\left(\mathcal{T}_-^+\right)^T |\mathbf{a}_{x_I}^+ \;=\; |\mathbf{a}_{x_I}^- \mathcal{T}_+^- \tag{58}$$

on interpolation operators $\mathcal{T}_+^-$ and $\mathcal{T}_-^+$. Transposing the relation in (58) leads to:

$$|\mathbf{a}_{x_I}^+ \mathcal{T}_-^+ \;=\; \left(\mathcal{T}_+^-\right)^T |\mathbf{a}_{x_I}^-. \tag{59}$$

We note here that similar constraints relating the interpolation operators from the two sides of the interface have been drawn in existing literature, see, for instance, equation (15) of [29] and equation (1) of [30].

Given the relations in (58) and (59), we observe that by choosing the set of parameters $\sigma_P^-$, $\sigma_V^-$, $\sigma_P^+$ and $\sigma_V^+$ such that linear system

$$\begin{cases} 1+\sigma_P^- + \sigma_V^- &= 0 \\ 1+\sigma_P^+ + \sigma_V^+ &= 0 \\ \sigma_V^+ - \sigma_P^- &= 0 \\ \sigma_P^+ - \sigma_V^- &= 0 \end{cases} \tag{60}$$

holds, we have ${}^{dE}/_{dt} = 0$ and therefore, semi-discretized system (55) is energy-conserving. For instance, $\sigma_P^- = \sigma_V^- = \sigma_P^+ = \sigma_V^+ = -{}^1/_2$ is one such set of parameters. Later, we indeed choose the parameters as such for the numerical examples presented in Section 4. Moreover, in Appendix B, we provide pairs of interpolation operators, i.e., $\mathcal{T}_+^-$ and $\mathcal{T}_-^+$, that provide at least second-order accurate interpolation and satisfy the relation in (58), for a collection of rational grid spacing ratios.

### *3.2.4. Heterogeneous media*

In this section, we allow the physical parameters $\rho$ and $c$ of the acoustic wave equation to be heterogeneous and demonstrate that the results obtained in previous sections can be applied in this case straightforwardly. Specifically, we consider the following PDE system:

$$\begin{cases} \frac{1}{\rho c^2}\frac{\partial p}{\partial t} = -\left(\frac{\partial u}{\partial x} + \frac{\partial v}{\partial y}\right); & \text{(61a)} \\ \rho\frac{\partial u}{\partial t} = -\frac{\partial p}{\partial x}; & \text{(61b)} \\ \rho\frac{\partial v}{\partial t} = -\frac{\partial p}{\partial y}. & \text{(61c)} \end{cases}$$

Comparing to (1), the physical parameters are moved to the left-hand side of (61). As shown later, it is natural to perform energy analysis on this form. The potential and kinetic energy associated with (61) are given as

$$\mathscr{E}_p = \frac{1}{2}\int_{\square}\frac{1}{\rho c^2}p^2 dxdy; \tag{62a}$$

$$\mathscr{E}_k = \frac{1}{2}\int_{\square}\rho\left(u^2 + v^2\right)dxdy, \tag{62b}$$

respectively. The total energy associated with (61) is the sum of $\mathscr{E}_p$ and $\mathscr{E}_k$:

$$\mathscr{E} = \mathscr{E}_p + \mathscr{E}_k. \tag{63}$$

We start with discretization on a staggered uniform grid, cf. Figure 1, and assume that these physical parameters are discretized on the same subgrids that their corresponding solution variables are associated with, i.e., $1/\rho c^2$ in (61a) is discretized on the $p$ subgrid, $\rho$ in (61b) is discretized on the $u$ subgrid, and $\rho$ in (61c) is discretized on the $v$ subgrid. Using the discretization operators derived in Section 3.2.1, we obtain the following semi-discretized system:

$$\begin{cases} C^P\mathcal{A}^P\frac{dP}{dt} = -\mathcal{A}^P\mathcal{D}_x^U U - \mathcal{A}^P\mathcal{D}_y^V V; & \text{(64a)} \\ C^U\mathcal{A}^U\frac{dU}{dt} = -\mathcal{A}^U\mathcal{D}_x^P P; & \text{(64b)} \\ C^V\mathcal{A}^V\frac{dV}{dt} = -\mathcal{A}^V\mathcal{D}_y^P P, & \text{(64c)} \end{cases}$$

where $C^P$, $C^U$ and $C^V$ are diagonal matrices having the same sizes as $\mathcal{A}^P$, $\mathcal{A}^U$ and $\mathcal{A}^V$, respectively, standing for discretization of $1/\rho c^2$ on the $p$ subgrid, discretization of $\rho$ on the $u$ subgrid and discretization of $\rho$ on the $v$ subgrid, respectively. We refer to these matrices as coefficient matrices hereafter. Since both $C^P$ and $\mathcal{A}^P$ are diagonal, their product $C^P\mathcal{A}^P$ is also diagonal and hence symmetric, which can be used to define the discrete counterpart of the potential energy, cf. (62a). Similarly, $C^U\mathcal{A}^U$ and $C^V\mathcal{A}^V$ can be used to define the discrete counterpart of the kinetic energy, cf. (62b).

We define the discrete energy associated with (64) as

$$E = \frac{1}{2}P^T\left(C^P\mathcal{A}^P\right)P + \frac{1}{2}U^T\left(C^U\mathcal{A}^U\right)U + \frac{1}{2}V^T\left(C^V\mathcal{A}^V\right)V, \tag{65}$$

which mimics the continuous energy given in (63). Taking the temporal derivative of (65) and substituting in (64), we obtain exactly the same result as obtained in (44) and those follow it. Therefore, we conclude that the SBP operators derived in Section 3.2.1 and the SATs derived in Section 3.2.2 can be directly applied to system (61), leading to energy-conserving discretization of the acoustic wave equation with heterogeneous media.

Similarly, for discretization on staggered block-wise uniform grid with nonconforming interface, cf. Figure 6, the procedure and results of Section 3.2.3 also apply directly to (61), after appropriate modifications to the definition of discrete energy.

## 4. Numerical examples

In this section, we present several numerical examples to validate the techniques derived in previous sections by solving PDE system (1) with different sets of physical parameters and discretization configurations. The SBP operators developed in Section 3.2.1 are used to discretize the spatial derivatives in (1). These SBP operators are then modified by SATs, as demonstrated in Sections 3.2.2 and 3.2.3, to account for boundary and interface conditions, whenever necessary. We denote the resulting spatial discretization operators as $\tilde{\mathcal{D}}_x^U$, $\tilde{\mathcal{D}}_y^V$, $\tilde{\mathcal{D}}_x^P$ and $\tilde{\mathcal{D}}_y^P$, which correspond to derivatives $\partial u/\partial x$, $\partial v/\partial y$, $\partial p/\partial x$ and $\partial p/\partial y$, respectively.

After applying the aforementioned spatial discretization, we use the staggered leapfrog scheme to integrate over time, as summarized in Algorithm 1, where the superscripts attached to $P$, $U$ and $V$ indicate the time steps. Variable $p$ is discretized at integer time steps while variables $u$ and $v$ are discretized at midpoints between integer time steps.

**Algorithm 1** Staggered Leapfrog Scheme

**for** $i_t = 1, \ldots, n_t$ **do**

*At half step, update the velocities:*

$$U^{(i_t+1/2)} = U^{(i_t-1/2)} + \Delta t \left(C^U\right)^{-1} \tilde{\mathcal{D}}_x^P P^{(i_t)};$$

$$V^{(i_t+1/2)} = V^{(i_t-1/2)} + \Delta t \left(C^V\right)^{-1} \tilde{\mathcal{D}}_y^P P^{(i_t)};$$

*At full step, update the pressure:*

$$P^{(i_t+1)} = P^{(i_t)} + \Delta t \left(C^P\right)^{-1} \left(\tilde{\mathcal{D}}_x^U U^{(i_t+1/2)} + \tilde{\mathcal{D}}_y^V V^{(i_t+1/2)}\right) + \Delta t\, \mathcal{S}^{(i_t+1/2)};$$

**end for**

We have presented techniques to fully discretize PDE system (1) and are now ready to address the actual physical and numerical configurations.

### 4.1. Block-wise homogeneous media

In this example, we consider a two-layer medium that is piecewise homogeneous. The medium parameters within the top layer are given by $\rho = 0.5\text{kg/m}^3$ and $c = 1\text{m/s}$. The medium parameters within the bottom layer are given by $\rho = 1\text{kg/m}^3$ and $c = 2\text{m/s}$. We use a point source $\mathcal{S}$ to drive the wave propagation, whose temporal profile is the Ricker wavelet with central frequency of 5Hz and time delay of 0.25s. We count the maximal frequency of the source content as 12.5Hz.

For spatial discretization, we consider a finite difference grid consisting of two horizontal layers, maintaining a 2:1 ratio of the grid spacing. On the vertical direction, the lowest grid row of the top layer and the highest grid row of the bottom layer overlap. Grid spacings are decided on a points-per-wavelength basis. In both layers, we ask to have 10 grid points per minimal wavelength, which translates to grid spacings $\Delta_x^F = 0.008\text{m}$ for the top layer and $\Delta_x^C = 0.016\text{m}$ for the bottom layer. The top layer is composed of a $p$ subgrid of size $120 \times 61$, a right shifted $u$ subgrid of size $120 \times 61$ and a up shifted $v$ subgrid of size $120 \times 60$. Similarly, the bottom layer is composed of a $p$ subgrid of size $60 \times 31$, a right shifted $u$ subgrid of size $60 \times 31$ and an up shifted $v$ subgrid of size $60 \times 30$. The grid configuration is illustrated in Figure 6 for a reduced number of grid points.

We complete PDE system (1) with periodic boundary condition on left and right boundaries, and free surface boundary condition on top and bottom boundaries. Moreover, the two grid layers are coupled together by the interface condition (53). We use fourth-order staggered central differences to approximate the spatial derivatives in the interior of the simulation domain.

In $x$-direction, when approaching the left and right boundaries, the stencils are wrapped around so that the periodic boundary condition is accounted for. We note here that by design (cf. Figure 6), the grids start with $p$ subgrid points and end with $u$ subgrid points in $x$-direction so that periodic boundary condition can be imposed most naturally in the staggered setting. With this treatment, $\mathbf{q}_x$ of (35) is a zero matrix, so does $Q_x$ since $Q_x = \mathbf{q}_x \otimes \underline{\mathbf{a}}_y$, cf. (43b). Consequently, no term will pop out in the 2D discrete energy analysis due to the periodic boundary condition. In $y$-direction, when approaching boundaries or interface, the stencils are adapted as in (12a) and (12b). The interface condition and free surface boundary condition are imposed via the SATs developed in Section 3.2.3 and Section 3.2.2, respectively. The interpolation formulas used on the interface are presented in (B.1) and (B.2) of Appendix B.

Both the source and receiver locations are placed on the $p$ subgrid. Specifically, we place the point source at $5\Delta_x^F$ rightward from the left boundary and $5\Delta_x^F$ downward from the top boundary to drive the wave propagation, and place the receiver at $5\Delta_x^F$ leftward from the right boundary and $5\Delta_x^F$ downward from the top boundary to record the time history of pressure response. The obtained signal is referred to as the seismogram in the following.

With the above specifications, we simulate the wave propagation for 5000 time steps with time step length 0.0012s, i.e., 6s in total. The obtained seismogram, referred to as the SBP result in the following, is displayed in Figure 7 and compared against finite element simulation result. Here, we compare against finite element simulation result, instead of uniform grid finite difference simulation result, because of the jumps in medium parameters. From Figure 7, we observe excellent agreement between the SBP result and the FEM result.

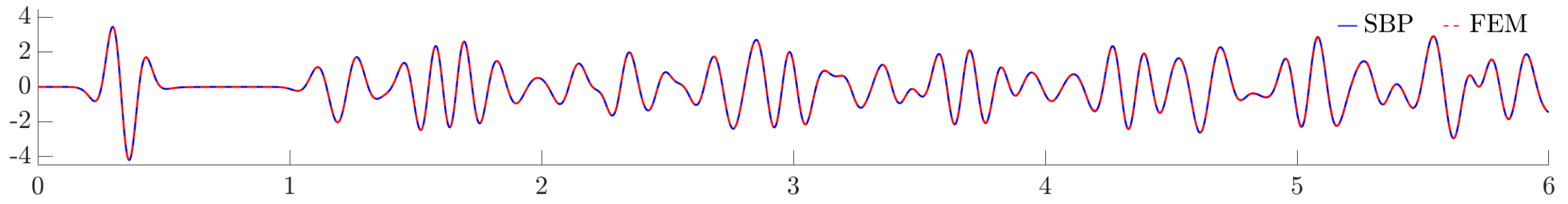

Figure 7: Seismograms simulated for the block-wise homogeneous model with the SBP and FEM discretization for 6s.

Now, we move to examining the long time stability of the SBP discretization, which is of special concern for this work. We run the simulation with the SBP discretization for 60s and display the pressure response of the last 6s in Figure 8. We observe stable behavior with no visible trace of unstable modes as those appearing in Figure 3. In Figure 9, the discrete energy given in (65) is displayed.[4] We observe that after the source effect tapers off, the discrete energy remains steady, confirming the energy-conserving property developed throughout Section 3.2.

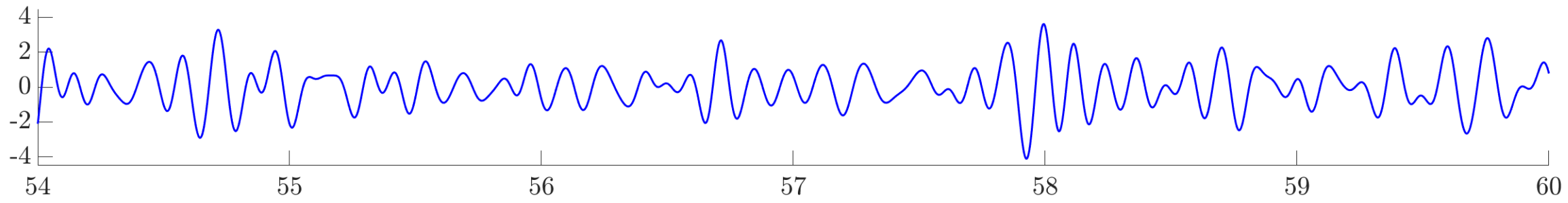

Figure 8: The last 6s of the 60s seismogram simulated for the block-wise homogeneous model with the SBP discretization.

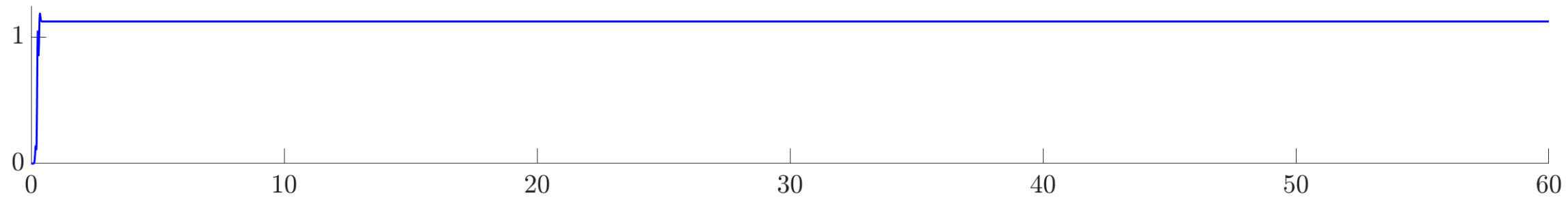

Figure 9: Discrete energy associated with the block-wise homogeneous model and the SBP discretization for 60s.

### 4.2. *Laterally homogeneous media*

In this example, we consider a laterally homogeneous medium that varies smoothly with depth. Specifically, the density $\rho$ increases from 0.5kg/m$^3$ (top boundary) to 1kg/m$^3$ (bottom boundary) linearly and the wave speed $c$ increases from 1m/s (top boundary) to 2m/s (bottom boundary) linearly. We consider a finite difference grid consisting of two horizontal layers, maintaining a 6:5 ratio of the grid spacing. The minimal and maximal wave speeds within the top layer are 1m/s and 1.2m/s, respectively. The minimal and maximal wave speeds within the bottom layer are 1.2m/s and 2m/s, respectively. Grid spacings are decided on a points-per-wavelength basis. In both layers, we ask to have 10 grid points per minimal wavelength, which translates to grid spacings $\Delta_x^F = 0.008$m for the top layer and $\Delta_x^C = 0.0096$m for the bottom layer.

The top layer is composed of a $p$ subgrid of size $120 \times 25$, a right shifted $u$ subgrid of size $120 \times 25$ and an up shifted $v$ subgrid of size $120 \times 24$. The bottom layer is composed of a $p$ subgrid of size $100 \times 81$, a right shifted $u$ subgrid of size $100 \times 81$ and an up shifted $v$ subgrid of size $100 \times 80$. The interpolation formulas used on the interface are presented in (B.6) of Appendix B. The rest of the numerical specifications (e.g., boundary and interface conditions, source and receiver, time integration) are identical to those used in Section 4.1.

We simulate the wave propagation for 6s and display the recorded pressure response in Figure 10. This result, referred to as the nonuniform result, is compared against the uniform grid finite difference simulation result of the same problem, which is also displayed in Figure 10. As in the previous example, we observe excellent agreement between the two seismograms.

[4]Due to grid staggering in time, the variable $p$ is not discretized at the same time instances as those for variables $u$ and $v$. As remedy, we average values of the discretized vectors of $p$ from two neighboring time steps before evaluating the discrete energy according to (65).

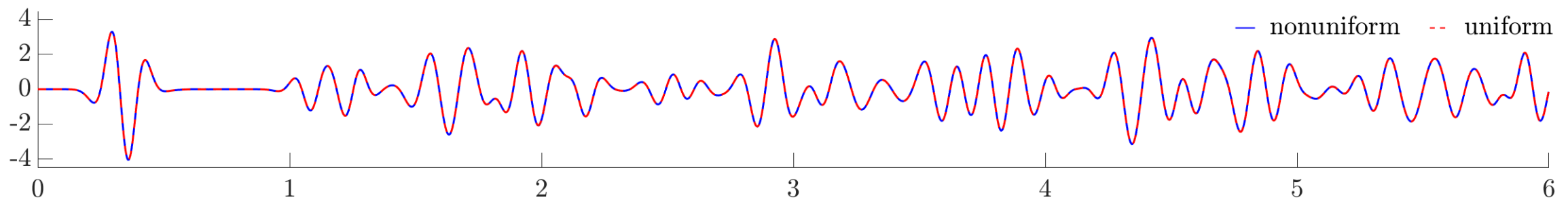

Figure 10: Seismograms simulated for the laterally homogeneous model with the nonuniform and uniform finite difference discretization for 6s.

In Figures 11-12, we demonstrate the long time stable behavior of the nonuniform grid simulation. Figure 11 displays last 6s of the pressure response simulated with the nonuniform grid discretization for 60s. Figure 12 displays the discrete energy, defined by (65), of the laterally homogeneous model. As in the previous example, these figures demonstrate the stable and energy-conserving properties developed throughout Section 3.2.

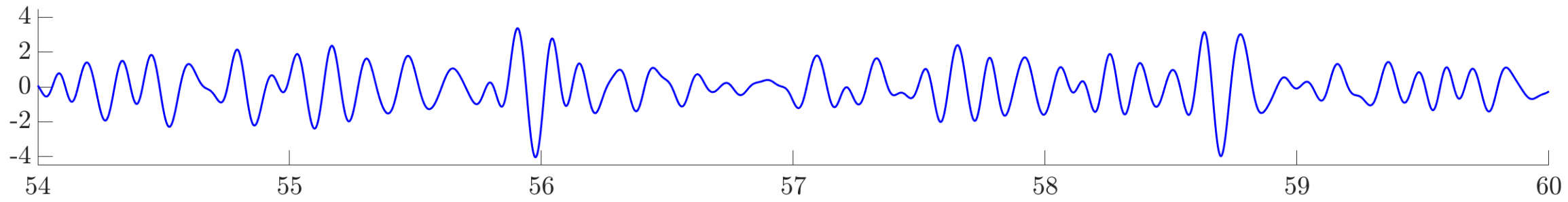

Figure 11: Last 6s of the 60s seismogram simulated for the laterally homogeneous model with the nonuniform grid discretization

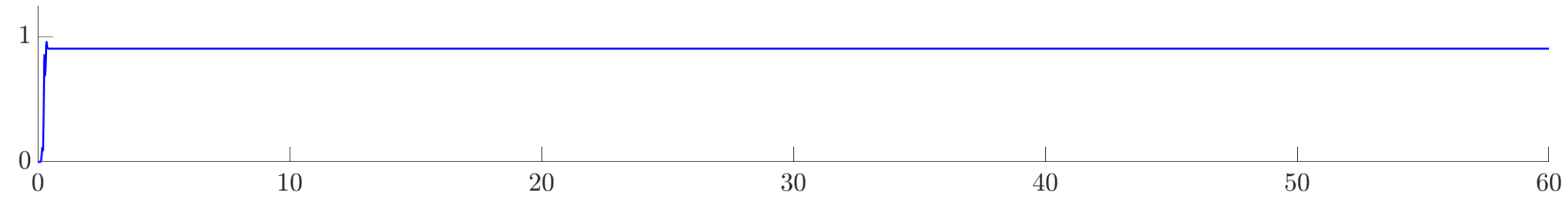

Figure 12: Discrete energy associated with the laterally homogeneous model and the nonuniform grid discretization for 60s.

### 4.3. *The Marmousi model*

In this example, we consider the Marmousi model, cf. [31], which is a standard test case for seismic wave modeling. Profiles of the density ($\rho$) and wave speed ($c$) of the Marmousi model are demonstrated in Figure 13. Both parameters are characterized by 751 rows and 2301 columns of data points, uniformly distributed with 4m distance between two neighboring rows and columns.

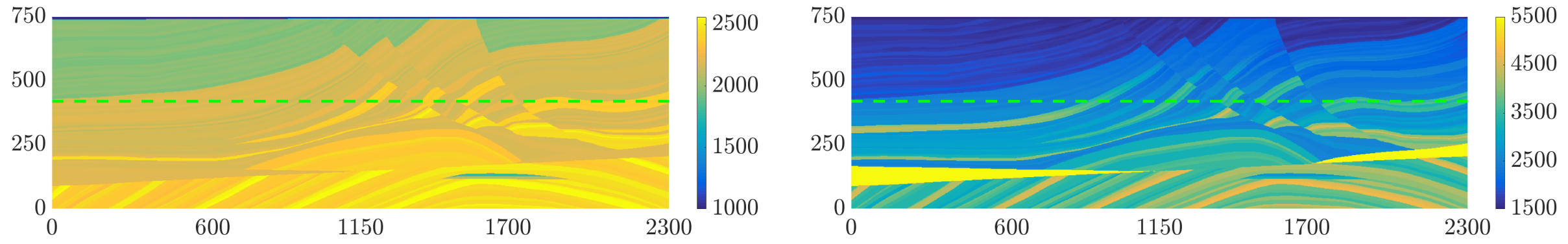

(a) Density (kg/m$^3$).

(b) Velocity (m/s).

Figure 13: Parameters of the Marmousi model.

We consider a finite difference grid consisting of two horizontal layers, maintaining a 3:2 ratio of the grid spacing. The top layer is composed of a $p$ subgrid of size $2301 \times 331$, a right shifted $u$ subgrid of size $2301 \times 331$ and an up shifted $v$ subgrid of size $2301 \times 330$. The bottom layer is composed of a $p$ subgrid of size $1534 \times 281$, a right shifted $u$ subgrid of size $1534 \times 281$ and an up shifted $v$ subgrid of size $1534 \times 280$.

The grid spacing used for the top layer is $\Delta_x^F = 4$m, matching the distance between the data points of the Marmousi model. The grid spacing used for the bottom layer is $\Delta_x^C = 6$m. When a finite difference grid point does not match a data point in location, we interpolate the model parameters with bilinear interpolation.

We use a point source to drive the wave propagation, whose temporal profile is the Ricker wavelet with central frequency of 12Hz and time delay of 0.25s. Maximal frequency of the source content is counted as 30Hz. Boundary and interface conditions are identical to those used in the previous examples and are dealt with in the same manner. The interpolation formulas used on the interface are presented in (B.3) of Appendix B.

In Table 1, discretization information for both layers are listed, including minimal wave speed $v_{min}$, maximal wave speed $v_{max}$, their ratio, minimal wavelength $\lambda_{min}$ (corresponding to the maximal frequency 30Hz), grid spacing $\Delta x$ and grid points per minimal wavelength, in that order. As comparison, the same information for a uniform grid finite difference discretization is also listed therein.

| | $v_{min}$ | $v_{max}$ | $v_{min}/v_{max}$ | $\lambda_{min}$ | $\Delta x$ | ppmw |
|---|---|---|---|---|---|---|
| Top | 1500m/s | 3650m/s | (≈) 0.411 | 50m | 4m | 12.5 |
| Bottom | 2405m/s | 5500m/s | (≈) 0.437 | (≈) 80.17m | 6m | (≈) 13.36 |
| Uniform | 1500m/s | 5500m/s | (≈) 0.273 | 50m | 4m | 12.5 (top)<br>(≈) 20.04 (bottom) |

Table 1: Discretization information for the Marmousi model. The term 'ppmw' in the last column is the abbreviation of 'points per minimal wavelength'.

Comparing to the uniform grid discretization, the benefit of the nonuniform grid discretization is twofold. First, the spatial sampling rates, i.e., grid points per minimal wavelength, for the two layers are kept close to each other, cf. last column of Table 1. In contrast, for the uniform grid case, the number of grid points per minimal wavelength (20.04) is significantly larger than that of the top layer region (12.5), implying that the bottom layer region is oversampled in space. Second, the grid spacing $\Delta x$ is usually determined by specifying the number of grid points per minimal wavelength $\lambda_{min}$. Using $N_{ppmw}$ to denote this number, we can write $\Delta x = \lambda_{min}/N_{ppmw} = v_{min}/(f_{max} \cdot N_{ppmw})$. Moreover, the time step length $\Delta t$ is usually restricted by the Courant-Friedrichs-Lewy (CFL) condition, i.e,

$$\Delta t \leq C \frac{\Delta x}{v_{max}} = \frac{C}{f_{max} \cdot N_{ppmw}} \frac{v_{min}}{v_{max}},$$

where $C$ is a constant associated with the particular problem. We observe that the CFL constraint on time step length is proportional to the ratio $v_{min}/v_{max}$. The smaller the ratio is, the stronger $\Delta t$ is constrained. As evidenced by the fourth column of Table 1, using nonuniform grid can alleviate the CFL constraint.

In Figure 14, pressure responses of 3s are displayed for both the nonuniform grid and the uniform grid simulations. The time step length is chosen as 2.5e-4s for both cases, for the convenience of comparison. As in the previous examples, we observe excellent agreement of the two seismograms. Moreover, Figure 15 displays the last 3s of the pressure response from the nonuniform grid simulation for 60s with time step length 3.75e-4s. Figure 16 displays the discrete energy, defined by (65), of the Marmousi model. As in the previous examples, these figures demonstrate the stable and energy-conserving properties developed throughout Section 3.2.

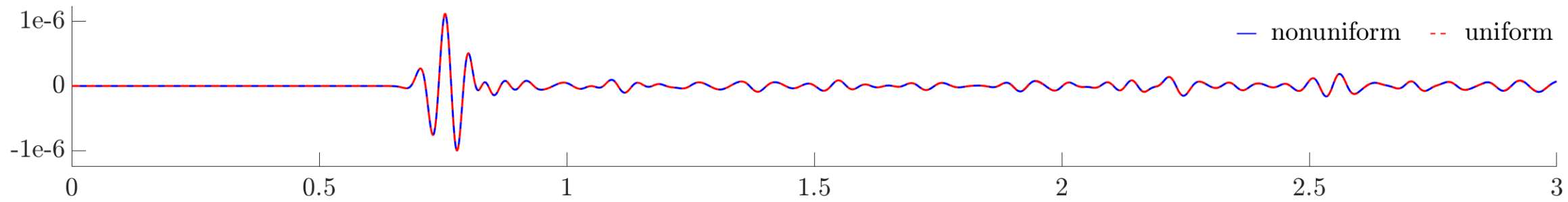

Figure 14: Seismograms simulated for the Marmousi model with the nonuniform and uniform finite difference discretization for 3s.

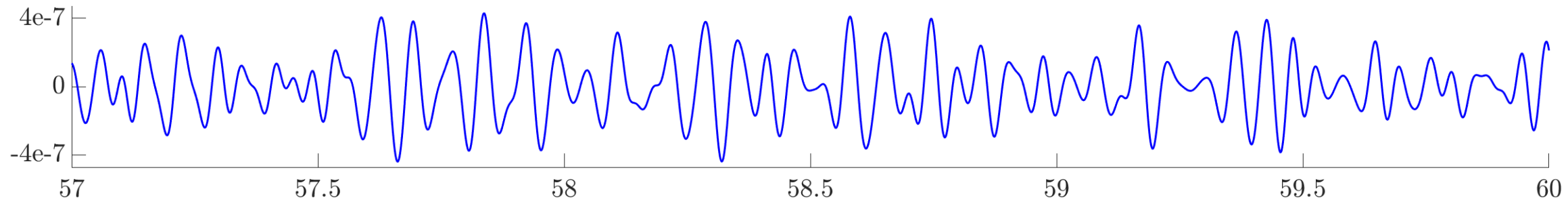

Figure 15: The last 3s of the 60s seismogram simulated for the Marmousi model with the nonuniform grid discretization.

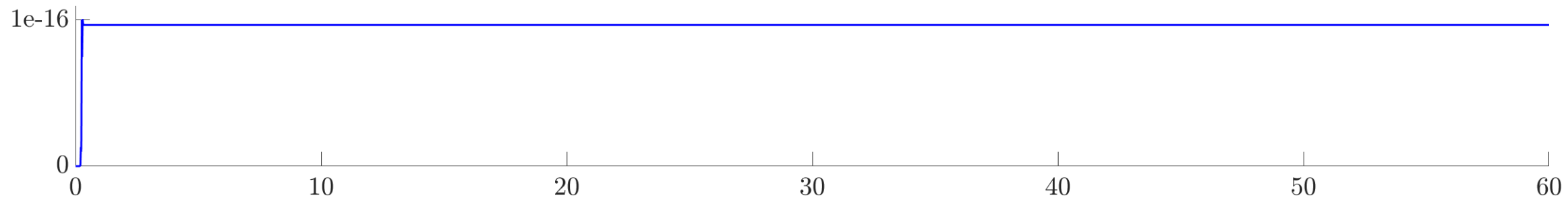

Figure 16: Discrete energy associated with the Marmousi model and the nonuniform grid discretization for 60s.

## 5. Several remarks

In this section, we make several remarks on issues that are not the focus of this work, but are of practical interests for seismic wave propagation, or other applications in general.

### *5.1. CFL constraint*

The SBP-SAT discretization can impose stricter constraint on time step length than that imposed by the underlying stencil used in interior of the simulation domain. For instance, considering the 1D wave equation (2) discretized by the staggered fourth-order central difference stencil $[{}^{1}/_{24},\ {}^{-9}/_{8},\ {}^{9}/_{8},\ {}^{-1}/_{24}]/\Delta x$ and the staggered leapfrog time integration scheme, it can be shown (see, e.g., [9]) that the time step length is limited by $\Delta t \leq \frac{6}{7}\Delta x \approx 0.8571\,\Delta x$. If the same PDE system, accompanied by the free surface boundary condition on both end-points, is discretized by the SBP operators in (12), the SATs in Section 3.1.2 and the staggered leapfrog time integration scheme, numerical experiments reveal that the maximal time step length allowed is roughly $0.635\,\Delta x$. Similar observation can be made for the 2D case, where the numbers become 0.6061 versus 0.5105. For this particular case, ratio of the grid spacing does not seem to affect the limit. In general, the limit on time step length remains an underexplored topic for the SBP-SAT discretization. Interested readers may consult [32], where free parameters left in the design of SBP operators are used to minimize the spectral radii of the SBP operators so that larger time step length can be used in simulation.

### *5.2. Point source*

Since the SATs technique weakly imposes boundary condition by penalty, the boundary condition may not be strictly satisfied when strong variation in the solution appears nearby the boundary. The discrepancy from the boundary condition then acts as an error source and affects the accuracy in the interior. This is indeed the case when a point source is placed too close to the boundary. However, this error diminishes quickly as the mesh is refined. Therefore, one can devise local refinement technique to deal with this issue. Alternatively, one may consult the techniques presented in [33]. These techniques distribute the effect of the point source over a collection of grid points, which can be biased away from the boundary, to achieve high order accurate discretization for the point source.

### *5.3. Perfectly matched layer*

Perfectly matched layer (PML), cf. [34, 35], is a technique commonly used in the simulation of propagating waves of various types. PMLs absorb outgoing waves with little reflection so that an infinite or semi-infinite physical domain can be mimicked by a finite simulation domain. Interested readers may consult [36], which addresses PML within the SBP-SAT discretization framework.

### *5.4. Future prospects*

First, the techniques presented in this work can be extended to the 3D case straightforwardly. As demonstrated by the 2D case in Section 3.2, the building blocks for higher dimensional SBP discretization operators are the SBP discretization operators for the corresponding 1D systems, e.g., (31) and (34). In other words, spatial derivatives in different directions are separate concerns in the process of constructing SBP operators. Moreover, proper 2D interface transfer operators can be constructed via tensor product of their 1D counterparts presented in this work. Relation (58) will be carried over to these 2D operators due to the properties of tensor product.

Also, the presented techniques can be generalized and applied to elastic wave equations naturally. Specifically, the elastic constitutive relations can be written in terms of the compliance tensor, which effectively moves the material parameters to the left-hand sides of the equations, as in (61) for acoustic wave equations. Moreover, the compliance tensor enters the definition of the strain energy (potential energy) naturally, similar to the term ${}^{1}/_{\rho c^2}$ in (62a) for acoustic wave equations. Interested readers may consult [18, Appendix B] for more detail. With this maneuver, the workflow presented in this work can be mapped to the elastic case identically.

## 6. Conclusion

In this work, we study the numerical simulation of acoustic wave equations arising from seismic applications, with special focus on staggered grid finite difference methods. We explore mechanisms that can extend the finite difference methods to block-wise uniform grids in a stable manner so that wave speeds in the subterranean media, which can have large variations, can be efficiently accounted for. We use the concepts of summation-by-parts (SBP) and simultaneous approximation terms (SATs) to achieve this. Specifically, we design finite difference discretization operators satisfying the SBP property for spatial derivatives appearing in the acoustic wave equation. These operators are applied in each block independently. The neighboring blocks are coupled together through the SATs, which act as penalty terms and impose the interface conditions weakly. Specially designed interpolation formulas on the interface are integral parts of the SATs. The overall discretization strategy is shown to be accurate and energy-conserving, which is verified by numerical examples of both theoretical and practical values.

## 7. Acknowledgments

Gao and Keyes gratefully acknowledge the support of KAUST's Office of Sponsored Research under CCF-CAF/URF/1-2596.

## Appendix A. Convergence tests

In this section, we test the accuracy of the discretization operators defined in (37) through mesh refinement. Specifically, we consider the acoustic wave equation with unit constant coefficients, cf. (28), defined over the square domain $(0,1)\times(0,1)$. The following manufactured solution:

$$\begin{cases} p(x,y,t) &= \sin(4\pi x)\sin(4\pi y)\cos(4\pi\sqrt{2}t); \\ u(x,y,t) &= -\frac{\sqrt{2}}{2}\cos(4\pi x)\sin(4\pi y)\sin(4\pi\sqrt{2}t); \\ v(x,y,t) &= -\frac{\sqrt{2}}{2}\sin(4\pi x)\cos(4\pi y)\sin(4\pi\sqrt{2}t), \end{cases} \tag{A.1}$$

which is taken from [20], is used for the convergence tests here. We note here that since both $\sin(4\pi x)$ and $\cos(4\pi x)$ span two full periods over interval $[0,1]$, the aforementioned manufactured solution satisfies periodic boundary condition on the $x$-direction. Moreover, since $p(x,0,t) = p(x,1,t) = 0$, it also satisfies free surface boundary conditions on the top and bottom boundaries. These will indeed be the boundary conditions used in the upcoming tests.

First, we consider the discretization scenario of Section 3.2.2, where the square domain $(0,1)\times(0,1)$ is uniformly discretized. In the interior of the square domain, we use the fourth-order staggered central difference stencil $[1/24,\ -9/8,\ 9/8,\ -1/24]/\Delta x$ to approximate the spatial derivatives in (28). When approaching the left and right boundaries, the stencil is wrapped around to account for periodic boundary condition in $x$-direction. When approaching the top and bottom boundaries, the stencil is adapted as in (12a) and (12b), for $\partial v/\partial y$ and $\partial p/\partial y$, respectively. Diagonal components of the 1D norm matrices $|\mathbf{a}_x$ and $\mathbf{a}|_x$ are simply the grid spacing $\Delta x$ for all grid points. Diagonal components of the 1D norm matrices $\underline{\mathbf{a}}_y$ and $\overline{\mathbf{a}}_y$ are again $\Delta x$ for interior grid points, but adapt for grid points near the top and bottom boundaries as in (12c) and (12d). Free surface boundary conditions on the top and bottom boundaries are imposed via the SATs technique developed in Section 3.2.2.

We run simulations for (28) on successively refined grids of sizes $16\times16$, $32\times32$, $64\times64$ and $128\times128$ for 1 second with time step length $\Delta t$ = 1E-06 second. The second-order staggered leapfrog scheme is used for time integration. We measure the error at the final time step with the weighted $\ell_2$ norm, with the norm matrices $\mathcal{A}^P$, $\mathcal{A}^U$ and $\mathcal{A}^V$ being the weights. These errors and their associated convergence rates are displayed in Table A.1.

| | 16 × 16 | 32 × 32 | 64 × 64 | 128 × 128 | 256 × 256 |
|---|---|---|---|---|---|
| Error | 4.20E-02 | 4.13E-03 | 3.39E-04 | 2.90E-05 | 2.53E-06 |
| Conv. Rate | — | 3.34 | 3.61 | 3.55 | 3.51 |

Table A.1: Errors and convergence rates for uniform discretization.

Next, we consider the discretization scenario of Section 3.2.3, where the square domain $(0, 1) \times (0, 1)$ is split into two blocks separated by a horizontal interface placed at $y = 0.5$. Both blocks are discretized uniformly, but with different grid spacings, relating to each other with a 2:1 ratio. The discretization operators for each block are built in the same way as in the previous case. So do the boundary treatment. The interface conditions, cf. (53), are imposed via the SATs technique developed in Section 3.2.3. Specifically, the interface interpolation operators are given in (B.1) and (B.2) of Appendix B.

Again, we run simulations for (28) on successively refined grids, which consist of $16 \times 8$, $32 \times 16$, $64 \times 32$ and $128 \times 64$ grid points for the bottom block and, correspondingly, $32 \times 16$, $64 \times 32$, $128 \times 64$ and $256 \times 128$ grid points for the top block. Parameters concerning time integration are the same as those of the previous case. Errors measured in the weighted $\ell_2$ norm and their associated convergence rates are displayed in Table A.2.

| | $16 \times 8$<br>$32 \times 16$ | $32 \times 16$<br>$64 \times 32$ | $64 \times 32$<br>$128 \times 64$ | $128 \times 64$<br>$256 \times 128$ | $256 \times 128$<br>$512 \times 256$ |
|---|---|---|---|---|---|
| Error | 3.89E-02 | 3.84E-03 | 3.15E-04 | 2.70E-05 | 2.37E-06 |
| Conv. Rate | — | 3.34 | 3.61 | 3.54 | 3.51 |

Table A.2: Errors and convergence rates for block-wise uniform discretization.

Since derivative approximations at grid points near the top and bottom boundaries, as well as the interface, are second-order accurate by design (cf. Section 3.1.1), third-order accuracy can be expected for the overall discretization (cf. [37, 38]), which is confirmed by the convergence rates in Tables A.1 and A.2. Moreover, we observe that the convergence rates are not degraded by the interface interpolation operators.

## Appendix B. Formulas for interface transfer operators

In this section, we present interface transfer operators that satisfy the accuracy constraints and the additional constraint (58) for various ratios of the grid spacing. We note here that interpolation formulas satisfying these constraints are not unique. The formulas we present retain symmetry from the grid configuration, whenever available, and in general have short lengths.

Since in $x$-direction, the grids considered in this article are uniform and subjected to periodic boundary condition, we assume that both $|\mathbf{a}^+_{x_I}$ and $|\mathbf{a}^-_{x_I}$ are identity matrices scaled by their respective grid spacing. Moreover, we assume that grid points from both sides of the interface match at certain locations. We define the elemental interval as the smallest interval whose two endpoints are locations where grid points from the two sides of the interface match. For instance, in Figures B.1–B.6, the intervals between the two dashed lines are elemental intervals. Since the grids, and hence the interpolation formulas, repeat themselves on different elemental intervals, we only present the interpolation formulas for these elemental intervals in the following.

We start with the case of 2:1 ratio of grid spacing and consider two interpolation problems depicted in Figures B.1 and B.2, respectively. For the problem illustrated in Figure B.1, we aim to approximate the values of smooth function $\phi$ at fine grid points $x^F_0$, $x^F_1$ and $x^F_2$ (blue solid points in Figure B.1), given its values at coarse grid points $x^C_{-1}$, $x^C_0$, $x^C_1$ and $x^C_2$ (red solid points in Figure B.1). On the other hand, for the reciprocal problem illustrated in Figure B.2, we aim to approximate the values of $\phi$ at coarse grid points $x^C_0$ and $x^C_1$ (red solid points in Figure B.2), given its values at fine grid points (blue solid points in Figure B.1). These two interpolation problems are linked together by the additional constraint given in (58).

Formulas in (B.1) and (B.2) provide third-order accurate interpolation results for these two interpolation problems, respectively. Moreover, the interface transfer operators induced from (B.1) and (B.2) satisfy the additional constraint (58). In (B.1), column vector $\Phi^C = \left[\phi(x^C_{-1}),\ \phi(x^C_0),\ \phi(x^C_1),\ \phi(x^C_2)\right]^T$ is composed of values of $\phi$ at the coarse grid points with which we interpolate. On the other hand, column vector $\tilde{\Phi}^F = \left[\tilde{\phi}(x^F_0),\ \tilde{\phi}(x^F_1),\ \tilde{\phi}(x^F_2)\right]^T$ is composed of approximations of $\phi$ at fine grid points, i.e., the interpolation results. Column vectors $\Phi^F$ and $\tilde{\Phi}^C$ in (B.2) are defined similarly. These formulas are derived with the assistance of the symbolic computing software Maple.

$$\tilde{\Phi}^F = \begin{bmatrix} 0 & 1 & 0 & 0 \\ -1/16 & 9/16 & 9/16 & -1/16 \\ 0 & 0 & 1 & 0 \end{bmatrix} \Phi^C. \tag{B.1}$$

$$\tilde{\Phi}^C = \begin{bmatrix} -1/32 & 0 & 9/32 & 1/2 & 9/32 & 0 & -1/32 & 0 & 0 \\ 0 & 0 & -1/32 & 0 & 9/32 & 1/2 & 9/32 & 0 & -1/32 \end{bmatrix} \Phi^F. \tag{B.2}$$

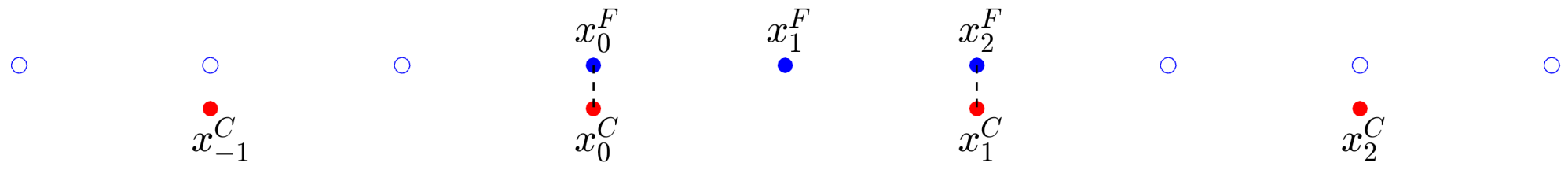

Figure B.1: Grid configuration associated with formulas (B.1).

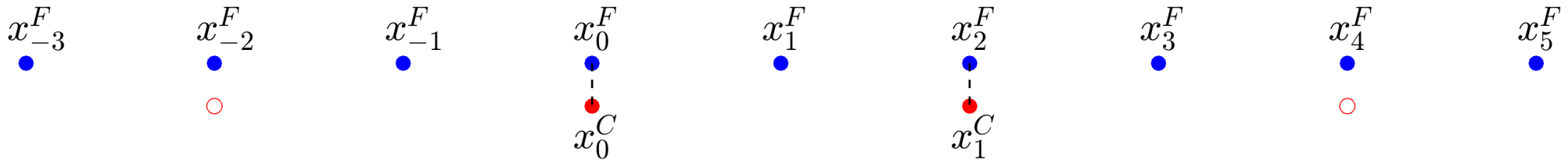

Figure B.2: Grid configuration associated with formulas (B.2).

We note here that since the relation in (58) is built into these formulas, the two sets of interpolation formulas can be recovered from each other. Therefore, for other ratios of grid spacing addressed in the following, we only present the interpolation formulas from coarse grid points to fine grid points in order to save space. Specifically, formulas in (B.3)-(B.6) provide at least second-order accurate interpolation results for the interpolation problems depicted in Figures B.3-B.6, respectively, addressing ratios of grid spacing 3:2, 4:3, 5:4 and 6:5, respectively. Their counterparts also provide at least second-order accurate interpolation results for the reciprocal interpolation problems of those depicted in Figures B.3-B.6, respectively.

$$\tilde{\Phi}^F = \begin{bmatrix} -1/96 & 1/24 & 15/16 & 1/24 & -1/96 & 0 & 0 \\ 0 & -13/288 & 103/288 & 217/288 & -19/288 & 0 & 0 \\ 0 & 0 & -19/288 & 217/288 & 103/288 & -13/288 & 0 \\ 0 & 0 & -1/96 & 1/24 & 15/16 & 1/24 & -1/96 \end{bmatrix} \Phi^C. \tag{B.3}$$

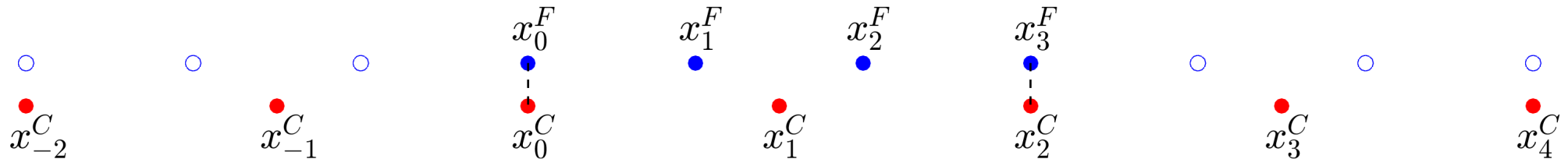

Figure B.3: Grid configuration associated with formulas (B.3).

$$\tilde{\Phi}^F = \begin{bmatrix} 0 & 0 & 1 & 0 & 0 & 0 & 0 & 0 \\ -1/288 & -13/864 & 2/9 & 83/96 & -59/864 & 0 & 0 & 0 \\ 0 & -1/432 & -1/18 & 241/432 & 241/432 & -1/18 & -1/432 & 0 \\ 0 & 0 & 0 & -59/864 & 83/96 & 2/9 & -13/864 & -1/288 \\ 0 & 0 & 0 & 0 & 0 & 1 & 0 & 0 \end{bmatrix} \Phi^C. \tag{B.4}$$

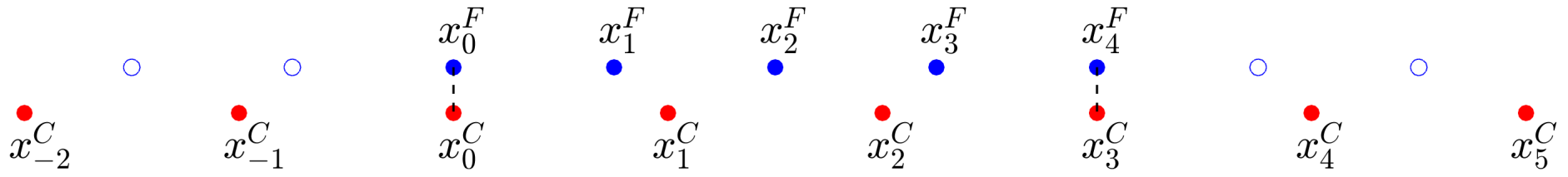

Figure B.4: Grid configuration associated with formulas (B.4).

$$\tilde{\Phi}^F = \begin{bmatrix} -63/3200 & 63/800 & 1411/1600 & 63/800 & -63/3200 & 0 & 0 & 0 & 0 \\ 0 & -123/3200 & 777/3200 & 2631/3200 & -61/3200 & -3/400 & 0 & 0 & 0 \\ 0 & 0 & -183/3200 & 29/64 & 531/800 & -93/1600 & -1/640 & 0 & 0 \\ 0 & 0 & -1/640 & -93/1600 & 531/800 & 29/64 & -183/3200 & 0 & 0 \\ 0 & 0 & 0 & -3/400 & -61/3200 & 2631/3200 & 777/3200 & -123/3200 & 0 \\ 0 & 0 & 0 & 0 & -63/3200 & 63/800 & 1411/1600 & 63/800 & -63/3200 \end{bmatrix} \Phi^C. \tag{B.5}$$

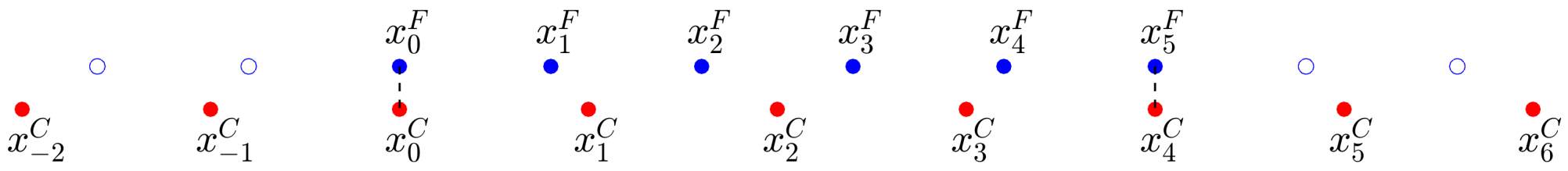

Figure B.5: Grid configuration associated with formulas (B.5).

$$\tilde{\Phi}^F = \begin{bmatrix} -39/2000 & 39/500 & 883/1000 & 39/500 & -39/2000 & 0 & 0 & 0 & 0 & 0 \\ 0 & -631/18000 & 317/1500 & 3781/4500 & -17/2250 & -161/18000 & 0 & 0 & 0 & 0 \\ 0 & 1/3600 & -317/6000 & 6871/18000 & 13043/18000 & -23/450 & -1/375 & 0 & 0 & 0 \\ 0 & 0 & 0 & -1/16 & 9/16 & 9/16 & -1/16 & 0 & 0 & 0 \\ 0 & 0 & 0 & -1/375 & -23/450 & 13043/18000 & 6871/18000 & -317/6000 & 1/3600 & 0 \\ 0 & 0 & 0 & 0 & -161/18000 & -17/2250 & 3781/4500 & 317/1500 & -631/18000 & 0 \\ 0 & 0 & 0 & 0 & 0 & -39/2000 & 39/500 & 883/1000 & 39/500 & -39/2000 \end{bmatrix} \Phi^C. \tag{B.6}$$

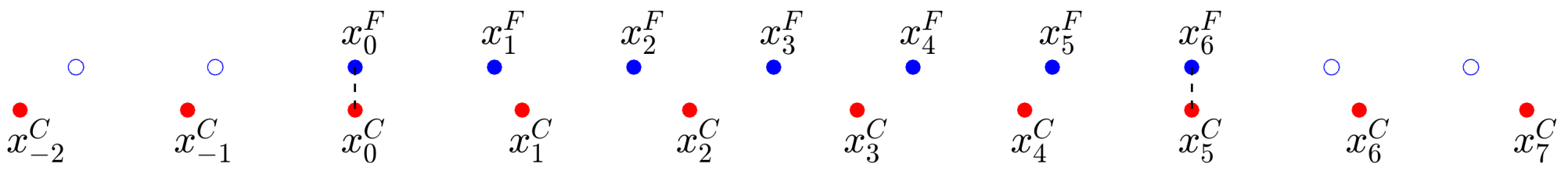

Figure B.6: Grid configuration associated with formulas (B.6).